\documentclass[11pt,a4paper]{article}
\usepackage{amssymb,amsmath,amsthm,epsfig,verbatim,pstricks,pst-node,url}
\usepackage{tikz,pgfplots,tikz-3dplot,psfrag}
\usepackage{graphbox} 
\usepackage{enumitem} 
\usepackage{mhchem} 
\usepackage{ifpdf}  
\newtheorem{theorem}{\sc Theorem.}[section]
\newtheorem{lemma}[theorem]{\sc Lemma.}
\newtheorem{remark}[theorem]{\sc Remark.}

\newenvironment{AMS}%
{{\upshape\bfseries AMS subject classifications. }\ignorespaces}{}
\newenvironment{keywords}{{\upshape\bfseries Key words. }\ignorespaces}{}

\newcommand{\bRplus}{{\mathbb R}_{>0}}
\newcommand{\bRgeq}{{\mathbb R}_{\geq 0}}

\newcommand{\tanspace}{{\rm T}}        
\newcommand{\RZ}{{\mathbb R} \slash {\mathbb Z}}
\newcommand{\bR}{{\mathbb R}}

\newcommand{\matpartnpar}[1]{\partial_{#1}^{\text{\tiny$\square$}}}
\newcommand{\matpartn}{\matpartnpar{t}}

\newcommand{\erf}{\operatorname{erf}}

\newcommand{\sign}{\operatorname{sign}}

\newcommand{\dH}[1]{\;{\rm d}{\mathcal{H}}^{#1}} 

\newcommand{\id}{{\rm id}}

\newcommand{\dd}[1]{\frac{\rm d}{{\rm d}#1}}
\newcommand{\ddt}{\dd{t}}

\newcommand{\ek}{e}
\newcommand{\Skappa}{H}
\def\epsilon{\varepsilon} 

\def\hat{\widehat}

\newcommand{\errorSS}{\|\mathcal{S} - \mathcal{S}_h\|_{L^\infty}}

\newcommand{\errorXX}{\|x - x_h\|_{L^\infty}}
\begin{document}
\title{
On the numerical approximation of \\ hyperbolic mean curvature flows for
surfaces}
\author{Klaus Deckelnick\footnotemark[2]\ \and 
        Robert N\"urnberg\footnotemark[3]}

\renewcommand{\thefootnote}{\fnsymbol{footnote}}
\footnotetext[2]{Institut f\"ur Analysis und Numerik,
Otto-von-Guericke-Universit\"at Magdeburg, 39106 Magdeburg, Germany \\
{\tt klaus.deckelnick@ovgu.de}}
\footnotetext[3]{Dipartimento di Mathematica, Universit\`a di Trento,
38123 Trento, Italy \\ {\tt robert.nurnberg@unitn.it}}

\date{}

\maketitle

\begin{abstract}
The paper addresses  the numerical approximation of two variants of hyperbolic mean curvature flow of surfaces in $\mathbb R^3$. For each evolution law
we propose both a finite element method, as well as a finite difference scheme in the case of axially symmetric surfaces. We present a number of numerical
simulations, including convergence tests as well as simulations suggesting the onset of singularities.
\end{abstract} 


\renewcommand{\thefootnote}{\arabic{footnote}}

\setcounter{equation}{0}
\section{Introduction} 

In this paper, we are going to investigate the numerical
approximation of hyperbolic mean curvature flows of the form
\begin{equation} \label{eq:hypgen}
\matpartn \mathcal{V}_{\mathcal S} = g(\mathcal V_{\mathcal S}^2) \Skappa  \quad \text{ on } \mathcal S(t)
\end{equation}
for a family of closed oriented surfaces $(\mathcal{S}(t))_{t\in[0,T]}$ in $\bR^3$.
Here $\mathcal{V}_{\mathcal S}$ denotes the velocity of $(\mathcal{S}(t))_{t\in[0,T]}$
in the direction of the unit normal $\nu_{\mathcal S}$,
$\matpartn$ is the normal time derivative on $(\mathcal{S}(t))_{t\in[0,T]}$,
and $\Skappa$ denotes the mean curvature of $\mathcal{S}(t)$ with the sign convention that the unit sphere with outward normal satisfies
$\Skappa = -2$. 
As for the coefficient function $g: \bRgeq \to \bRplus$,
we shall focus on two choices: Firstly,
the choice $g(s) = 1$ leads to the equation
\begin{equation} \label{eq:Gurtinbeta0}
\matpartn \mathcal{V}_{\mathcal S} = \Skappa \quad\text{on }\ \mathcal{S}(t),
\end{equation}
which is a special case of an evolution law derived by Gurtin and Podio--Guidugli in \cite{GurtinP91} as a model for the
evolution of melting-freezing waves at the solid-liquid interface of crystals such as \ce{^{4}He} helium. The equation \eqref{eq:Gurtinbeta0}
corresponds to \cite[(1.7)]{GurtinP91} for the case of an isotropic surface energy, unit effective density, vanishing kinetic coefficient  and  absence of external forcings.   \\
Secondly, the choice $g(s)=1+\tfrac12 s$ leads to the evolution law
\begin{equation} \label{eq:LeFloch}
\matpartn \mathcal{V}_{\mathcal{S}} = \tfrac12 
\left( \mathcal{V}_{\mathcal S}^2 + 2 \right)
\Skappa \quad \text{on }\ \mathcal{S}(t),
\end{equation}
which has been derived by LeFloch and Smoczyk in \cite{LeFlochS08} with the help of the
Hamiltonian principle. As a result, solutions of \eqref{eq:LeFloch} conserve the sum $E(t)=\tfrac12 \int_{\mathcal S(t)} (\mathcal{V}_{\mathcal S}^2 + 2) \dH2$ of kinetic and internal energy,  see Corollary 4.7 in \cite{LeFlochS08}. Here, ${\rm d}\mathcal{H}^2$ denotes integration with respect to the 
two-dimensional Hausdorff measure. Recalling \cite[Theorem~32]{bgnreview}, we can derive this property from \eqref{eq:LeFloch} via
\begin{align}  \label{eq:ELeFloch}
 E'(t) & = \tfrac12 \int_{\mathcal S(t)} \bigl( \matpartn ( \mathcal V_{\mathcal S}^2 +2) - (\mathcal V_{\mathcal S}^2 +2) \mathcal V_{\mathcal S} H \bigr) \dH2 \nonumber \\
& = \int_{\mathcal S(t)} \mathcal V_{\mathcal S} \matpartn \mathcal V_{\mathcal S} \dH2 - \tfrac12 \int_{\mathcal S(t)} ( \mathcal V_{\mathcal S}^2 + 2) \mathcal V_{\mathcal S} H \dH2 =0,
\end{align}
while on the other hand we see that solutions of \eqref{eq:Gurtinbeta0} satisfy
\begin{equation} \label{eq:EGurtinbeta0}
\ddt \int_{\mathcal S(t)} e^{\frac12\mathcal{V}_{\mathcal S}^2} \dH2 = \int_{\mathcal S(t)} e^{\frac12 \mathcal V_{\mathcal S}^2}  \mathcal V_{\mathcal S} \matpartn \mathcal V_{\mathcal S} \dH2 -   \int_{\mathcal S(t)} e^{\frac12 \mathcal V_{\mathcal S}^2}  \mathcal V_{\mathcal S} H \dH2 = 0.
\end{equation}
Thus, even though the two flows \eqref{eq:Gurtinbeta0} and \eqref{eq:LeFloch} are derived in different ways, they share similar conservation properties in that the integrand $\tfrac12 \mathcal{V}_{\mathcal S}^2 + 1$ in $E(t)$ may be 
interpreted as a low order approximation of the integrand 
$e^{\frac12\mathcal{V}_{\mathcal S}^2}$ in \eqref{eq:EGurtinbeta0}.

Analytical results for the flow \eqref{eq:Gurtinbeta0} have focussed so far on the evolution of planar curves, see \cite{KongLW09,KongW09}. It is shown in \cite{KongLW09}
that if the initial curve is strictly convex and the initial velocity vector does not point outwards, then the evolution exists on a finite time interval $[0,T_{max})$. Furthermore,
as $t \nearrow T_{max}$, the solution either shrinks to a point or converges to a convex curve with discontinuous curvature. Local well-posedness for the flow \eqref{eq:LeFloch}
has been obtained in \cite{LeFlochS08} along with criteria that guarantee a blow-up in finite time. 
Let us also refer to \cite{Yau00} and \cite{Notz13} for 
other hyperbolic systems that
are related to the subject of our paper but are not considered here.
\\
Let us next address the numerical approximation of hyperbolic mean curvature flows. A finite difference scheme in order to approximate planar curves evolving by  \eqref{eq:Gurtinbeta0},
which is based on a parametric approach, has been proposed and analyzed in \cite{hmcf}. In the semidiscrete case, optimal error estimates are obtained as long as the
solution is smooth. Numerical experiments suggest the occurrence of singularities in finite time thereby numerically confirming the analytical results of \cite{KongLW09}. The authors in  \cite{RotsteinBN99} 
develop an algorithm for the evolution of planar, polygonal curves according to a version of \eqref{eq:Gurtinbeta0} in the case of a crystalline, anisotropic surface energy. 
In addition, in \cite{GinderS16}
an implicit description of curves evolving by a variant of \eqref{eq:Gurtinbeta0} is used, and an algorithm based on thresholding for the linear wave equation
is presented. 
Let us finally mention that geometric second order hyperbolic PDEs have 
recently been used in \cite{DongHZ21} for applications in image processing. 

The aim of this paper is to develop numerical methods for the evolution of {\it surfaces} either by \eqref{eq:Gurtinbeta0} or by \eqref{eq:LeFloch}. Our primary approach is based on a
parametric description of the moving surfaces and discretizes a system of PDEs satisfied by the position vector of the surface. The corresponding systems for both evolution laws are introduced in 
Section~\ref{sec:mathform}, see Lemma~\ref{lem:xtnormal}. 
The systems are chosen in such a way that the
resulting flow is normal, i.e.\ the velocity vector has no tangential component. This property has proven to be helpful both in analytical investigations, see \cite{KongLW09,LeFlochS08},
and in the numerical analysis carried out in \cite{hmcf}. 
We remark that for related parabolic geometric evolution equations,
like mean curvature flow, surface diffusion and Willmore flow, several
different approaches on introducing suitable tangential motions to help
maintain or improve the mesh quality exist, see e.g.\
\cite{gflows3d,willmore,bgnreview,ElliottF17,HuL22,DuanL24}. However, for
simplicity we will focus on purely normal formulations for the hyperbolic flows
considered in this paper.
Our numerical schemes are introduced in Section~\ref{sec:disc}. We shall
adapt a finite element method due to Dziuk, \cite{Dziuk91}, for mean curvature flow, which approximates the evolving surface by a triangulated surface and parameterizes the solution in the next time step over the
current surface. Apart from this general case, we are also interested in the more special situation that the evolving surfaces are rotationally symmetric. Here,  the evolution can be described 
in terms of the profile curve and we are able to extend the finite difference approach developed in \cite{hmcf} to this setting. There are two reasons for considering the case of
axisymmetric surfaces: on the one hand, an algorithm adapted to this setting will substantially reduce the computational cost, while on the other hand we can use 
the corresponding numerical results in order to verify our finite element simulations for the general case. 
Those simulations will be presented in Section~\ref{sec:nr} for the two flows
\eqref{eq:Gurtinbeta0} and \eqref{eq:LeFloch}, covering both the general and the axisymmetric case. They show interesting oscillatory behaviour of the surfaces and
provide numerical evidence for the onset of singularities in three space dimensions.

\setcounter{equation}{0}
\section{Mathematical formulation} \label{sec:mathform}

\subsection{The general case}

In what follows we shall employ a parametric description of the
surfaces $\mathcal{S}(t)$ evolving by \eqref{eq:hypgen}, i.e.\
$\mathcal{S}(t) = \mathfrak{X}(\mathcal M,t)$ for some closed
reference surface $\mathcal M \subset \mathbb R^3$ and some mapping
$\mathfrak{X}:\mathcal M \times [0,T] \rightarrow \mathbb R^3$. We assume that
$\mathfrak{X}(\cdot,t)$ is a bijection between $\mathcal M$ and $\mathcal{S}(t)$
for each $t \in [0,T]$. 

Let us consider the following system 
\begin{subequations} \label{eq:xtt}
\begin{align}
\mathfrak{X}_{tt} \circ \mathfrak{X}^{-1} & = g(|\mathfrak{X}_t|^2\circ \mathfrak{X}^{-1}) \Skappa \nu_{\mathcal S} - \tfrac12 \nabla_{\mathcal{S}(t)} (|\mathfrak{X}_t|^2\circ \mathfrak{X}^{-1})
\qquad \text{on } \mathcal{S}(t) = \mathfrak{X}(\mathcal M,t), \label{eq:xtta} \\
\mathfrak{X}(\cdot,0) &= \mathfrak{X}_0, \quad \mathfrak{X}_t(\cdot,0)= V_0 \,  \nu_{\mathcal S} \circ \mathfrak{X}_0 \qquad \text{on } \mathcal M. \label{eq:xttb}
\end{align}
\end{subequations}
Here, $\mathfrak{X}_0:\mathcal M \rightarrow \mathbb R^3$ is a 
parameterization of the initial surface $\mathcal{S}_0$ and
$V_0: \mathcal M \rightarrow \mathbb R$ is a given initial normal velocity. Furthermore, $\nabla_{\mathcal S(t)} f$ denotes the tangential gradient of a differentiable
function $f:\mathcal S(t) \rightarrow \mathbb R$. 

\begin{lemma} \label{lem:xtnormal}
Let $g \in C^0(\bRgeq)$ and 
let $\mathfrak{X}: \mathcal{M} \times [0,T] \to \bR^3$ be a solution to \eqref{eq:xtt}. Then $\mathcal{S}(t)=
\mathfrak{X}(\mathcal M,t)$ satisfies \eqref{eq:hypgen}. Furthermore, $\mathfrak{X}$ is a normal 
parameterization, i.e.
\begin{equation} \label{eq:xtnormal}
\mathfrak{X}_t(p,t) \cdot \tau = 0 \qquad \forall\ \tau \in \tanspace_{\mathfrak{X}(p,t)}\,\mathcal{S}(t), \ \forall\ (p,t) \in \mathcal{M} \times [0,T],
\end{equation}
where $\tanspace_z \mathcal S(t)$ denotes the tangent space to $\mathcal S(t)$ at the point 
$z\in\mathcal{S}(t)$.
\end{lemma}
\begin{proof}
Let us fix $p \in \mathcal M$ and denote by $\phi:U \rightarrow \mathcal M$,
$U \subset \bR^2$, a local parameterization of $\mathcal M$
with $\phi(u)=p$. 
Abbreviating $\tilde{\mathfrak{X}}(u,t):=\mathfrak{X}(\phi(u),t)$ and observing that $\lbrace \tilde{\mathfrak{X}}_{u_1}(u,t), \tilde{\mathfrak{X}}_{u_2}(u,t) \rbrace$ is a 
basis of $\tanspace_{\mathfrak{X}(p,t)}\,\mathcal{S}(t)$ we obtain \eqref{eq:xtnormal} if we can show that
\begin{equation} \label{eq:nflow}
\tilde{\mathfrak{X}}_t(u,t) \cdot \tilde{\mathfrak{X}}_{u_i}(u,t)=0 \qquad \mbox{ for } i=1,2 \mbox{ and }0 \leq t \leq T.
\end{equation}
Since $\tilde{\mathfrak{X}}_{tt}(u,t)=\mathfrak{X}_{tt}(\phi(u),t)$ we obtain with the help of \eqref{eq:xtta} and the fact that $\nu_{\mathcal S} \cdot  \tilde{\mathfrak{X}}_{u_i}=0$ that
\begin{align}
(\tilde{\mathfrak{X}}_t \cdot \tilde{\mathfrak{X}}_{u_i})_t 
& = \tilde{\mathfrak{X}}_{tt} \cdot \tilde{\mathfrak{X}}_{u_i} + \tilde{\mathfrak{X}}_t \cdot \tilde{\mathfrak{X}}_{u_i,t} 
= - \tfrac12 \bigl( \nabla_{\mathcal{S}(t)} (|\tilde{\mathfrak{X}}_t|^2 \circ \tilde{\mathfrak{X}}^{-1}) \bigr) \circ \tilde{\mathfrak{X}} \cdot \tilde{\mathfrak{X}}_{u_i}
+ \tilde{\mathfrak{X}}_t \cdot  \tilde{\mathfrak{X}}_{u_i,t}.
\label{eq:klaus0}
\end{align}
Let us denote by $g_{ij}= \tilde{\mathfrak{X}}_{u_i} \cdot \tilde{\mathfrak{X}}_{u_j}$ the coefficients of the first fundamental form and write $(g^{ij})=(g_{ij})^{-1}$. Then 
\begin{align*}
\tfrac12 \bigl( \nabla_{\mathcal{S}(t)}( |\tilde{\mathfrak{X}}_t|^2 \circ \tilde{\mathfrak{X}}^{-1} ) \bigr) \circ \tilde{\mathfrak{X}} \cdot \tilde{\mathfrak{X}}_{u_i} &
 =  \tfrac12 \sum_{k,l=1}^2 g^{kl} (|\tilde{\mathfrak{X}}_t|^2)_{u_k} \tilde{\mathfrak{X}}_{u_l} \cdot \tilde{\mathfrak{X}}_{u_i} \nonumber \\
&= \sum_{k,l=1}^2 g^{kl} (\tilde{\mathfrak{X}}_t \cdot \tilde{\mathfrak{X}}_{t,u_k}) g_{li} 
 =  \sum_{k=1}^2 \delta_{ki}(\tilde{\mathfrak{X}}_t \cdot \tilde{\mathfrak{X}}_{t,u_k})
= \tilde{\mathfrak{X}}_t \cdot \tilde{\mathfrak{X}}_{u_i,t}
\end{align*}
which combined with \eqref{eq:klaus0} implies that $(\tilde{\mathfrak{X}}_t \cdot \tilde{\mathfrak{X}}_{u_i})_t =0$. We deduce with the help of \eqref{eq:xttb} that 
\begin{align}
\tilde{\mathfrak{X}}_t(u,0) \cdot \tilde{\mathfrak{X}}_{u_i}(u,0)= \mathfrak{X}_t(p,0) \cdot \tilde{\mathfrak{X}}_{u_i}(u,0) = V_0(p) \nu_{\mathcal S}(\mathfrak X_0(p)) \cdot \tilde{\mathfrak{X}}_{u_i}(u,0)=0,
\end{align}
which yields  \eqref{eq:nflow}. Next, we recall that the normal time derivative at $(x_0,t_0)$ with $x_0 \in \mathcal{S}(t_0)$ is given by 
$\matpartn f(x_0,t_0)=\ddt f(\gamma(t),t)_{|t=t_0}$, where 
$\gamma$ is the curve satisfying $\gamma'(t)=\mathcal V_{\mathcal{S}}(\gamma(t),t) \nu_{\mathcal S}(\gamma(t),t)$ and
$\gamma(t_0)=x_0$. Choosing $p \in \mathcal M$ with $\mathfrak{X}(p,t_0)=x_0$ we have that $\gamma(t)=\mathfrak{X}(p,t)$ since \eqref{eq:xtnormal} implies that
\begin{displaymath}
\gamma'(t)=\mathfrak{X}_t(p,t) = (\mathfrak{X}_t(p,t) \cdot (\nu_{\mathcal S} \circ \mathfrak{X})(p,t) ) (\nu_{\mathcal S} \circ \mathfrak{X})(p,t) = \mathcal V_\mathcal{S}(\gamma(t),t) \nu_{\mathcal S}(\gamma(t),t).
\end{displaymath}
Thus,
\[
\matpartn \mathcal{V}_{\mathcal{S}} \circ \mathfrak{X} = \ddt \mathcal V_{\mathcal{S}}(\gamma(t),t)= \bigl( \mathfrak{X}_t \cdot (\nu_{\mathcal S} \circ \mathfrak{X}) \bigr)_t = 
\mathfrak{X}_{tt} \cdot (\nu_{\mathcal S} \circ \mathfrak{X})  + \mathfrak{X}_{t}  \cdot (\nu_{\mathcal S} \circ \mathfrak{X})_t
= \mathfrak{X}_{tt}  \cdot  (\nu_{\mathcal S} \circ \mathfrak{X})  
\]
since $(\nu_{\mathcal S} \circ \mathfrak{X})_t(p,t)  \in T_{\mathfrak{X}(p,t)} \mathcal{S}(t)$. As a result we infer from \eqref{eq:xtta} that 
\begin{displaymath}
\matpartn \mathcal V_{\mathcal{S}} = \mathfrak{X}_{tt} \circ \mathfrak{X}^{-1} \cdot \nu_{\mathcal S} =  g(|\mathfrak{X}_t|^2\circ \mathfrak{X}^{-1})  \Skappa = g(\mathcal V_{\mathcal S}^2) \Skappa \qquad \mbox{ on } \mathcal{S}(t)
\end{displaymath}
since $| \mathfrak X_t |^2 = \bigl( \mathfrak X_t \cdot (\nu_{\mathcal S} \circ \mathfrak X) \bigr)^2 = \mathcal V_{\mathcal S}^2 \circ \mathfrak X$. Thus,  \eqref{eq:hypgen} holds.
\end{proof}

\subsection{The axisymmetric case}

\begin{figure}
\center
\newcommand{\AxisRotator}[1][rotate=0]{%
    \tikz [x=0.25cm,y=0.60cm,line width=.2ex,-stealth,#1] \draw (0,0) arc (-150:150:1 and 1);%
}
\begin{tikzpicture}[every plot/.append style={very thick}, scale = 1]
\begin{axis}[axis equal,axis line style=thick,axis lines=center, xtick style ={draw=none}, 
ytick style ={draw=none}, xticklabels = {}, 
yticklabels = {}, 
xmin=-0.2, xmax = 0.8, ymin = -0.4, ymax = 2.55]
after end axis/.code={  
   \node at (axis cs:0.0,1.5) {\AxisRotator[rotate=-90]};
   \draw[blue,->,line width=2pt] (axis cs:0,0) -- (axis cs:0.5,0);
   \draw[blue,->,line width=2pt] (axis cs:0,0) -- (axis cs:0,0.5);
   \node[blue] at (axis cs:0.5,-0.2){$\ek_1$};
   \node[blue] at (axis cs:-0.2,0.5){$\ek_2$};
   \draw[red,very thick] (axis cs: 0,0.7) arc[radius = 70, start angle= -90, end angle= 90];
   \node[red] at (axis cs:0.75,1.95){$\Gamma(t)$};
}
\end{axis}
\end{tikzpicture} \qquad \qquad
\tdplotsetmaincoords{120}{50}
\begin{tikzpicture}[scale=2, tdplot_main_coords,axis/.style={->},thick]
\draw[axis] (-1, 0, 0) -- (1, 0, 0);
\draw[axis] (0, -1, 0) -- (0, 1, 0);
\draw[axis] (0, 0, -0.2) -- (0, 0, 2.7);
\draw[blue,->,line width=2pt] (0,0,0) -- (0,0.5,0) node [below] {$\ek_1$};
\draw[blue,->,line width=2pt] (0,0,0) -- (0,0.0,0.5);
\draw[blue,->,line width=2pt] (0,0,0) -- (0.5,0.0,0);
\node[blue] at (0.2,0.4,0.1){$\ek_3$};
\node[blue] at (0,-0.2,0.3){$\ek_2$};
\node[red] at (0.75,0,1.9){$\mathcal{S}(t)$};
\node at (0.0,0.0,2.4) {\AxisRotator[rotate=-90]};

\tdplottransformmainscreen{0}{0}{1.4}
\shade[tdplot_screen_coords, ball color = red] (\tdplotresx,\tdplotresy) circle (0.7);
\end{tikzpicture}
\caption{Sketch of $\Gamma(t)$ and $\mathcal{S}(t)$, as well as 
the unit vectors $\ek_1$, $\ek_2$ and $\ek_3$.}
\label{fig:sketch}
\end{figure}
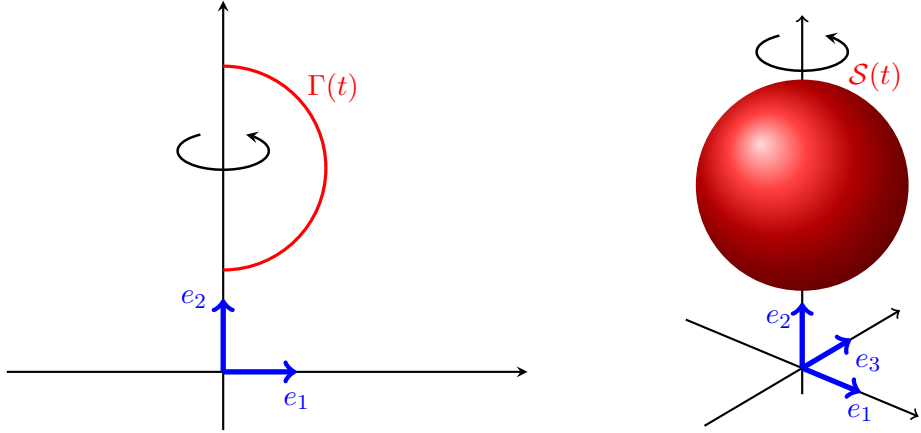

We would like to consider the above problems also in an axisymmetric setting. 
To do so, let 
\begin{displaymath}
I= \RZ,\ \partial I= \emptyset \quad \mbox{ or } \quad I=(0,1), \ \partial I= \lbrace 0,1 \rbrace.
\end{displaymath}
With a mapping $x:I \times [0,T] \rightarrow \mathbb R^2$ satisfying $x(\rho,t) \cdot e_1>0$ for $(\rho,t) \in I \times [0,T]$ we associate
the family of surfaces
\begin{equation} \label{eq:x2S}
\mathcal S(t)=\lbrace (x(\rho,t) \cdot e_1 \cos \theta, x(\rho,t) \cdot e_2, x(\rho,t) \cdot e_1 \sin \theta) : \rho \in I, \theta \in [0,2 \pi] \rbrace,
\end{equation}
which are generated by rotating the curves $\Gamma(t)=x(I,t)$ around the
$x_2$-axis, see also Figure~\ref{fig:sketch}, which is taken from 
\cite{axisd,aximcf}. 
In the case $I=\mathbb R/\mathbb Z$ the curve $\Gamma(t)$ is a closed curve lying in the half-space $\lbrace x \in \mathbb R^2 : x \cdot e_1>0 \rbrace$ 
and gives rise to a surface $\mathcal S(t)$ of genus 1. On the other hand, if $I=(0,1)$, we shall assume that 
\begin{displaymath}
x(\rho,t) \cdot e_1=0, \; x_\rho(\rho,t) \cdot e_2=0, \quad \rho \in \lbrace 0,1 \rbrace,
\end{displaymath}
so that the curve $\Gamma(t)$ meets the $x_2$-axis on the boundary vertically and $\mathcal S(t)$ is a surface of genus 0.
It can be shown that the normal velocity and the mean curvature of $\mathcal S(t)$ are given respectively by
\begin{equation} \label{eq:axiquant}
\mathcal V_{\mathcal S} = x_t \cdot \nu, \quad \Skappa = \kappa - \frac{\nu \cdot e_1}{x \cdot e_1},
\end{equation}
where $\nu=\tau^\perp$ is the unit normal and $\kappa$ is the curvature of $\Gamma(t)$. Here, $\tau= \frac{x_\rho}{| x_\rho |}$ is the unit tangent and 
$\cdot^\perp$ denotes the clockwise rotation through $\frac{\pi}{2}$.
See e.g.\ \cite{axisd,aximcf} for further details.

We propose to consider the following axisymmetric formulation for \eqref{eq:hypgen}:
\begin{subequations} \label{eq:app:hcsf}
\begin{alignat}{2}
& x_{tt} = g(| x_t|^2) \bigl(  \kappa \nu - \frac{\nu \cdot e_1}{x \cdot e_1} \nu \bigr)  -  (x_t \cdot \tau_t) \tau, \quad && \mbox{ in } I \times (0,T], \label{eq:app:hcsfa} \\
&  x \cdot e_1=0, \; x_\rho \cdot e_2=0, \qquad \quad && \mbox{ on } \partial I \times (0,T], \label{eq:app:hcsfb} \\[2mm]
 & x(\cdot,0) =x_0, \; x_t(\cdot,0) =V_0 \nu, \qquad && \mbox{ in } I. \label{eq:app:hcsfc}
 \end{alignat}
\end{subequations} 

Similarly to Lemma~\ref{lem:xtnormal} we have the following result:
\begin{lemma} \label{lem:xtnormalaxi}
Let $g \in C^0(\bRgeq)$ and 
let $x: I \times [0,T] \to \bR^2$ be a solution to \eqref{eq:app:hcsf}. Then 
the family of surfaces induced by \eqref{eq:x2S} 
satisfies \eqref{eq:hypgen}. Furthermore, $x$ is a normal parameterization, i.e.\ $x_t \cdot \tau = 0$ in $I \times [0,T]$.
\end{lemma}
\begin{proof} 
Let us begin by showing that $x$ is a normal parameterization. Using \eqref{eq:app:hcsfa} we obtain
\begin{displaymath}
(x_t \cdot \tau)_t = x_{tt} \cdot \tau + x_t \cdot \tau_t = - x_t \cdot \tau_t + x_t \cdot \tau_t=0,
\end{displaymath}
so that $x_t \cdot \tau =0$ since $(x_t \cdot \tau)_{|t=0} = (V_0 \nu) \cdot  \tau =0$ in view of \eqref{eq:app:hcsfc}. Next, in the same way 
as in Lemma~2.1 in \cite{hmcf} we infer using \eqref{eq:axiquant} together with the fact that $x_t \cdot \tau=0$ 
\begin{displaymath}
\matpartn \mathcal V_{\mathcal S} = (x_t \cdot \nu)_t = x_{tt} \cdot \nu = g(| x_t |^2) \bigl( \kappa - \frac{\nu \cdot e_1}{x \cdot e_1} \bigr)  = g(\mathcal V_{\mathcal S}^2) \Skappa,
\end{displaymath}
so that  \eqref{eq:hypgen} holds.
\end{proof}

\setcounter{equation}{0}
\section{Discretization} \label{sec:disc}

\subsection{The general case}

Let us begin by describing first the discretization with respect to time. 
To this end, we choose a time step size $\Delta t>0$ and set $t_m=m \Delta t$.
Let $\mathcal{S}^m \approx \mathcal{S}(t_m)$ as well as
$\mathfrak{X}^m \approx \mathfrak{X}(\cdot,t_m) : \mathcal{M} \to \bR^3$ with $\mathfrak{X}^m(\mathcal{M}) =
\mathcal{S}^m$. 
Following Dziuk's seminal paper on the numerical approximation 
for mean curvature flow, \cite{Dziuk91}, 
it will be convenient to work on the current
surface $\mathcal S^m$ rather than on the reference surface $\mathcal M$. 
Abbreviating 
$X^{m+1} = \mathfrak{X}^{m+1} \circ (\mathfrak{X}^m)^{-1}$ we have $X^{m+1} : \mathcal{S}^m \to \mathcal{S}^{m+1}$,
$X^{m} : \mathcal{S}^{m-1} \to \mathcal{S}^{m}$ as well as
\begin{displaymath}
\mathfrak{X}^{m-1} \circ (\mathfrak{X}^m)^{-1} = \bigl( \mathfrak{X}^m \circ ( \mathfrak{X}^{m-1})^{-1} \bigr)^{-1} = (X^m)^{-1},
\end{displaymath}
so that $(X^m)^{-1}$ parameterizes $\mathcal S^{m-1}$ over $\mathcal S^m$. 
The relationship between the different mappings is visualized below.
\ifpdf
\begin{center}
\includegraphics[angle=0,width=0.5\textwidth]{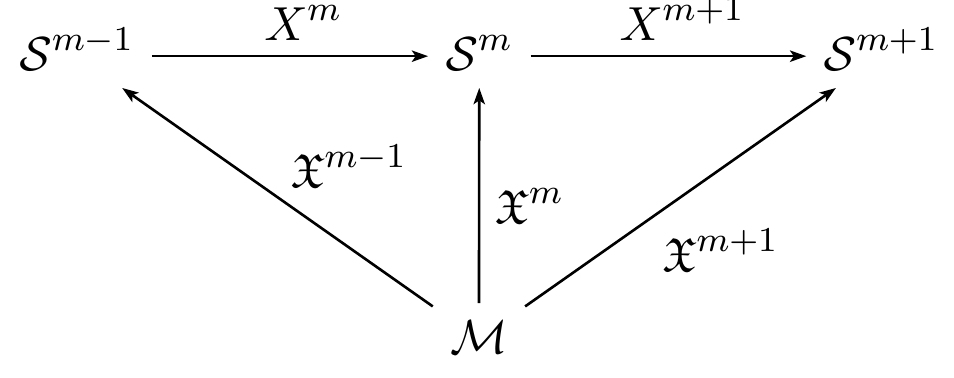}
\end{center}
\else
\[
\psset{arrows=->, arrowinset=0.25, linewidth=0.6pt, nodesep=4pt, labelsep=3pt, rowsep=1.8cm, colsep=2.5cm}
\begin{psmatrix}
  \mathcal{S}^{m-1} & \mathcal{S}^{m} & \mathcal{S}^{m+1} \\%
   & \mathcal{M} &
\end{psmatrix}
\ncline{1,1}{1,2}\naput[npos=0.55]{X^m}
\ncline{1,2}{1,3}\naput[npos=0.55]{X^{m+1}}
\ncline{2,2}{1,1}\nbput[npos=0.45]{\mathfrak{X}^{m-1}}
\ncline{2,2}{1,2}\nbput[npos=0.45]{\mathfrak{X}^m}
\ncline{2,2}{1,3}\nbput[npos=0.45]{\mathfrak{X}^{m+1}}
\]
\fi
It holds that 
\begin{equation} \label{eq:timediff2M}
\mathfrak{X}_{tt \mid_{t=t_m}} \approx \frac1{(\Delta t)^2} ( \mathfrak{X}^{m+1} - 2 \mathfrak{X}^m + \mathfrak{X}^{m-1})
\qquad \text{in } \mathcal M, 
\end{equation}
so that
\begin{align}
\mathfrak{X}_{tt \mid_{t=t_m}} \circ (\mathfrak{X}^m)^{-1} &\approx \frac1{(\Delta t)^2} 
( \mathfrak{X}^{m+1} \circ (\mathfrak{X}^m)^{-1} - 2\, \id + \mathfrak{X}^{m-1} \circ (\mathfrak{X}^m)^{-1}) \nonumber \\
& = 
\frac1{(\Delta t)^2} ( X^{m+1} - 2\, \id + (X^{m})^{-1} ) \qquad \text{on } \mathcal{S}^m.  \label{eq:timediff2}
\end{align}
Similarly, 
\begin{equation} \label{eq:timediff}
\mathfrak{X}_{t \mid_{t=t_m}} \circ (\mathfrak{X}^m)^{-1} \approx \frac1{\Delta t} 
( \id - (X^m)^{-1})
\qquad \text{on } \mathcal{S}^m.
\end{equation}
In order to write the term involving the mean curvature vector in divergence form, we observe that $\Skappa \nu_{\mathcal S} =\Delta_{\mathcal{S}(t_m)}\, \id$ 
on $\mathcal{S}(t_m)$, recall \cite{Dziuk91}, so that
\begin{align}
 g \bigl(  |\mathfrak{X}_{t}|^2 \circ \mathfrak{X}^{-1} \bigr) 
\Skappa \nu_{\mathcal S}  & = g \bigl( |\mathfrak{X}_{t}|^2 \circ \mathfrak{X}^{-1} \bigr) \nabla_{\mathcal S(t)} \cdot \nabla_{\mathcal S(t)} \id 
 \nonumber \\ 
& = \nabla_{\mathcal S(t)} \cdot \bigl( g( |\mathfrak{X}_{t}|^2 \circ \mathfrak{X}^{-1}) \nabla_{\mathcal S(t)} \id \bigr) - \nabla_{\mathcal S(t)} \id \,  \nabla_{\mathcal S(t)} g \bigl( |\mathfrak{X}_{t}|^2 \circ \mathfrak{X}^{-1} \bigr)
\nonumber \\
& = \nabla_{\mathcal S(t)} \cdot \bigl( g( |\mathfrak{X}_{t}|^2 \circ \mathfrak{X}^{-1}) \nabla_{\mathcal S(t)} \id \bigr) - g' \bigl(  |\mathfrak{X}_{t}|^2 \circ \mathfrak{X}^{-1} \bigr) \nabla_{\mathcal S(t)} \bigl( |\mathfrak{X}_{t}|^2 \circ \mathfrak{X}^{-1} \bigr),
\nonumber 
\end{align}
where we used that $(\nabla_{\mathcal S} \id)_{ij} =\delta_{ij} - \nu_{\mathcal S,i} \nu_{\mathcal S,j}$. Thus
\begin{align}  
& 
\bigl( g( |\mathfrak{X}_{t}|^2 \circ \mathfrak{X}^{-1}) \Skappa \nu_{\mathcal S}\bigr)_{\mid_{t=t_m}} \nonumber \\
& \approx 
\begin{cases}
\nabla_{\mathcal S^m} \cdot \nabla_{\mathcal S^m} \id & \mbox{if } g(s)=1, \\[2.5mm]
\tfrac12 \nabla_{\mathcal S^m} \cdot \bigl( ( |\mathfrak{X}_t|^2|_{t=t_m} \circ (\mathfrak{X^m})^{-1}+2) \nabla_{\mathcal S^m} \id \bigr) - \tfrac12 \nabla_{\mathcal S^m} \bigl( |\mathfrak{X}_t|^2_{|t=t_m} \circ (\mathfrak{X^m})^{-1} \bigr) & \mbox{if } g(s)=1+\tfrac12 s.
\end{cases}
\label{eq:timediff3}
\end{align}

Let us now combine \eqref{eq:timediff2}, \eqref{eq:timediff} and \eqref{eq:timediff3} to define time-discrete approximations for \eqref{eq:xtt}. In what follows it will be convenient to
replace $\nabla_{\mathcal S^m} \id$ by $\frac12\nabla_{\mathcal{S}^m} \bigl(X^{m+1} +  (X^m)^{-1}  \bigr)$ in \eqref{eq:timediff3} in order to enhance the stability and accuracy of the resulting fully discrete schemes. 
This approach is motivated by the authors' previous work in
\cite{hmcf}, as well as by the proposed numerical approximations 
for the linear wave equation in \cite[\S2.7]{GrossmannR07}. \\
Thus, we approximate \eqref{eq:xtta} in the case $g(s)=1$ by 
\begin{align}
&
\frac1{(\Delta t)^2} ( X^{m+1} - 2\, \id + (X^{m})^{-1} )
- \tfrac12 \nabla_{\mathcal{S}^m} \cdot \nabla_{\mathcal S^m} \bigl(X^{m+1} +  (X^m)^{-1}  \bigr) \nonumber \\ & \hspace{2cm}
= -\tfrac12 \nabla_{\mathcal{S}^m}  \left| \frac{\id - (X^{m})^{-1}}{\Delta t}\right|^2
\quad \text{on } \mathcal{S}^m, \label{eq:scheme1}
\end{align}
while in the case $g(s)=1+\tfrac12 s$ we choose the approximation
\begin{align}
&
\frac1{(\Delta t)^2} ( X^{m+1} - 2\, \id + (X^{m})^{-1} )
- \tfrac14 \nabla_{\mathcal{S}^m} \cdot  \left[  \Bigl( \left| \frac{\id - (X^{m})^{-1}}{\Delta t}\right|^2 +2 \Bigr) \nabla_{\mathcal S^m}  \bigl( X^{m+1} +  (X^m)^{-1} \bigr)  \right] \nonumber \\ & \hspace{2cm}
= - \nabla_{\mathcal{S}^m}  \left| \frac{\id - (X^{m})^{-1}}{\Delta t}\right|^2
\quad \text{on } \mathcal{S}^m. \label{eq:scheme2}
\end{align}
We now obtain the new approximation $\mathcal S^{m+1}$ of $\mathcal S(t_{m+1})$ from $\mathcal S^m$ and $(X^m)^{-1}$ by solving \eqref{eq:scheme1}, or \eqref{eq:scheme2}, for $X^{m+1}$
and then setting $\mathcal S^{m+1}=X^{m+1}(\mathcal S^m)$.  In order to start this time-discrete evolution,
$\mathcal S^0$ and $(X^0)^{-1}$ are required, which we obtain from the initial conditions \eqref{eq:xttb}.
Firstly, we let $\mathfrak X^0 = \mathfrak X_0$, which yields 
$\mathcal S^0=\mathfrak X^0(\mathcal M)$. In order to define
$\mathfrak X^{-1} \approx \mathfrak X(\cdot,-\Delta t)$, we
formally use a Taylor expansion to calculate
\begin{displaymath}
\mathfrak X(\cdot,-\Delta t) = \mathfrak X(\cdot,0) - \Delta t \mathfrak X_t(\cdot,0) + \tfrac12 (\Delta t)^2 \mathfrak X_{tt}(\cdot,0) + \mathcal{O}((\Delta t)^3),
\end{displaymath} 
and hence set 
$\mathfrak X^{-1} = \mathfrak X^0 - \Delta t \mathfrak X_t(\cdot,0) + \tfrac12 (\Delta t)^2 \mathfrak X_{tt}(\cdot,0)$. 
In view of \eqref{eq:xtt} we then have
\begin{align}
(X^0)^{-1} & = \mathfrak X^{-1} \circ (\mathfrak X^0)^{-1} 
= \id - \Delta t\, \mathfrak X_t(\cdot,0) \circ (\mathfrak X^0)^{-1} 
+ \tfrac12 (\Delta t)^2\, \mathfrak X_{tt}(\cdot,0) \circ
(\mathfrak X^0)^{-1} \nonumber \\
& = \id - \Delta t\, (V_0 \circ (\mathfrak X^0)^{-1}) \nu_{\mathcal S} 
+ \tfrac12 (\Delta t)^2 \, \bigl( g(V_0^2 \circ (\mathfrak X^0)^{-1}) \Skappa \nu_{\mathcal S} 
- \tfrac12 \nabla_{\mathcal S^0} (V_0 \circ (\mathfrak X^0)^{-1} )^2 \bigr) 
\ \mbox{on } \mathcal S^0. \nonumber
\end{align}
In the remainder of the paper, to simplify the presention, we will always assume that $V_0$ is constant so that the above relation simplifies to
\begin{equation} \label{eq:discinit}
(X^0)^{-1} = \id - \Delta t\, V_0 \nu_{\mathcal S} 
+ \tfrac12 (\Delta t)^2 \,  g(V_0^2) \Skappa \nu_{\mathcal S} 
\ \mbox{ on } \mathcal S^0. 
\end{equation}
Generalizing our schemes to a general initial velocity $V_0$ is
straightforward.

For a finite element discretization we closely follow \cite{Dziuk91}. 
Let $\mathcal{S}^{m}_h$ be a polyhedral surface,
i.e.\ a union of nondegenerate triangles with no hanging vertices
(see \cite[p.~164]{DeckelnickDE05}), approximating the
closed surface $\mathcal{S}(t_m)$.
Let $\mathcal{S}^m_h=\bigcup_{j=1}^J \overline{\sigma^m_j}$,
where $\{\sigma^m_j\}_{j=1}^J$ is a family of mutually disjoint open triangles
with vertices $\{p^m_k\}_{k=1}^K$ and set $h_m:=\max_{j=1\to J} {\rm diam}
(\sigma^m_j)$.
Then let
\begin{equation*}
\underline{V}(\mathcal{S}^m_h) := \{\chi \in C^0(\mathcal{S}^m_h,\bR^3): \chi
\!\mid_{\sigma^m_j} \mbox{ is affine}\ \forall\ j=1\to J\} 
=: [V(\mathcal{S}^m_h)]^3 ,
\label{eq:Vh}
\end{equation*}
where $V(\mathcal{S}^m_h)$ is the space of scalar continuous
piecewise linear functions on $\mathcal{S}^m_h$.
We define the natural Lagrange interpolation operator
$\pi^h_m : C^0(\mathcal{S}^m_h) \to V(\mathcal{S}^m_h)$ via
$(\pi^{m}_h\,\eta) (p^m_k) = \eta(p^m_k)$ for $k=1,\ldots,
K$. 
For functions $u,v\in L^2(\mathcal{S}^m_h)$,
we introduce the $L^2$ inner product 
$\langle\cdot,\cdot\rangle_{\mathcal{S}^m_h}$ via
\begin{equation*}
\langle u, v\rangle_{\mathcal{S}^m_h} = \int_{\mathcal{S}^m_h} u\,v \dH2.
\end{equation*}
In addition, if $u,\,v$ are piecewise continuous, with possible jumps
across the edges of $\{\sigma_j^m\}_{j=1}^J$,
we introduce the mass lumped inner product
$\langle\cdot,\cdot\rangle^h_{\mathcal{S}^m_h}$ as
\begin{equation}
\langle u, v \rangle^h_{\mathcal{S}^m_h} =
\tfrac13 \sum_{j=1}^J |\sigma^m_j|\sum_{k=0}^{2} 
(u\,v)((p^m_{j_k})^-),
\label{def:ip0}
\end{equation}
where $\{p^m_{j_k}\}_{k=0}^{2}$ are the vertices of 
$\sigma^m_j$, and where
we define $u((p^m_{j_k})^-) := \!
\underset{\sigma^m_j\ni p\to p^m_{j_k}}{\lim}\, u(p)$.
Moreover $|\sigma^m_j|$ denotes the measure of $\sigma^m_j$.
We generalize both inner products in the natural way to vector- and 
matrix-valued functions. In addition,
we introduce the outward unit normal $\nu^m_h$ to $\mathcal{S}^m_h$.

Our fully discrete approximation of \eqref{eq:xtta} is then defined as
follows.
Given $\mathcal{S}^m_h$ and 
$(X^{m}_h)^{-1} \in \underline{V}(\mathcal{S}^m_h)$, find
$X^{m+1}_h \in \underline{V}(\mathcal{S}^m_h)$ such that in the case $g(s)= 1$
\begin{align} \label{eq:fea_b}
&
\left\langle \frac{X^{m+1}_h - 2\, \id + (X^{m}_h)^{-1}}{(\Delta t)^2}
,\varphi_h \right\rangle^h_{\mathcal{S}^m_h}
+ \tfrac12\left\langle \nabla_{\mathcal{S}^m_h} (X^{m+1}_h + (X^{m}_h)^{-1})
,\nabla_{\mathcal{S}^m_h} \varphi_h \right\rangle_{\mathcal{S}^m_h} \nonumber \\ & \qquad
= -\tfrac12 \left\langle \nabla_{\mathcal{S}^m_h} 
\pi^m_h \left| \frac{\id - (X^{m}_h)^{-1}}{\Delta t}\right|^2,
\varphi_h \right\rangle_{\mathcal{S}^m_h}
\forall\ \varphi_h \in \underline{V}(\mathcal{S}^m_h),
\end{align} 
while in the case $g(s)=1+\tfrac12s$
\begin{align} \label{eq:feaLF}
&
\left\langle \frac{X^{m+1}_h - 2\, \id + (X^{m}_h)^{-1}}{(\Delta t)^2}
,\varphi_h \right\rangle^h_{\mathcal{S}^m_h} \nonumber \\ & \quad
+ \tfrac14\left\langle 
\left(\pi^m_h \left| \frac{\id - (X^{m}_h)^{-1}}{\Delta t}\right|^2 +2 \right)
\nabla_{\mathcal{S}^m_h} (X^{m+1}_h + (X^{m}_h)^{-1})
,\nabla_{\mathcal{S}^m_h} \varphi_h \right\rangle_{\mathcal{S}^m_h} \nonumber \\ & \qquad
= - \left\langle \nabla_{\mathcal{S}^m_h} 
\pi^m_h \left| \frac{\id - (X^{m}_h)^{-1}}{\Delta t}\right|^2,
\varphi_h \right\rangle_{\mathcal{S}^m_h}
\forall\ \varphi_h \in \underline{V}(\mathcal{S}^m_h).
\end{align} 
Having obtained $X^{m+1}_h$, we set $\mathcal S^{m+1}_h= X^{m+1}_h(\mathcal S^m_h)$.

\begin{lemma} \label{lem:ex}
Suppose that $\mathcal S^m_h$ consists of nondegenerate triangles. Then there exists a unique solution $X^{m+1}_h \in \underline{V}(\mathcal{S}^m_h)$
to \eqref{eq:fea_b} and \eqref{eq:feaLF}. 
\end{lemma}
\begin{proof} 
Since \eqref{eq:fea_b} is a linear scheme with the same number of unknowns as
equations, it is sufficient to prove that only the trivial solution solves
the homogeneous system
\[
\left\langle \frac{X_h}{(\Delta t)^2}
,\varphi_h \right\rangle^h_{\mathcal{S}^m_h}
+ \tfrac12\left\langle \nabla_{\mathcal{S}^m_h} X_h 
,\nabla_{\mathcal{S}^m_h} \varphi_h \right\rangle_{\mathcal{S}^m_h} 
= 0 \qquad \forall\ \varphi_h \in \underline{V}(\mathcal{S}^m_h) .
\]
Choosing $\varphi_h=X_h$ in the above immediately yields that $X_h=0$. The result for the scheme \eqref{eq:feaLF} follows in the same way.
\end{proof}

In order to define the initial data in the spirit of \eqref{eq:discinit}, we
first need to introduce continuous piecewise linear approximations of the
outer unit normal and the mean curvature vector. Following 
Definitions~51 and Equation~(29) in \cite{bgnreview}, 
we let $\omega^m_h \in \underline{V}(\mathcal{S}^m_h)$ be defined by
\begin{equation} \label{eq:omegah}
\left\langle \omega^m_h, \varphi_h \right\rangle_{\mathcal{S}^m_h}^h 
= \left\langle \nu^m_h, \varphi_h \right\rangle_{\mathcal{S}^m_h}
\quad\forall\ \varphi_h\in \underline{V}(\mathcal{S}^m_h),
\end{equation}
while $Y^m_h \in \underline{V}(\mathcal S^m_h)$ is the solution to
\begin{equation} \label{eq:kappavec}
\left\langle Y^m_h,\varphi_h \right\rangle^h_{\mathcal S^m_h}
+ \left\langle \nabla_{\mathcal S^m_h} \id,
\nabla_{\mathcal S^m_h} \varphi_h \right\rangle_{\mathcal S^m_h}
= 0 \quad \forall\ \varphi_h \in \underline{V}(\mathcal S^m_h).
\end{equation}
Now,  
let $\mathcal S^0_h$ be a triangular approximation of $\mathcal S^0$ and set
\begin{displaymath}
(X^0_h)^{-1}= \id - \Delta t\, V_0 \, \omega^0_h + \tfrac12 (\Delta t)^2\, g(V_0^2) 
Y^0_h \in \underline{V}(\mathcal S^0_h).
\end{displaymath}

\subsection{The axisymmetric case}

We shall employ a finite difference scheme in order to discretize 
\eqref{eq:app:hcsf} in space. To simplify the presentation, we only consider
the case $I=(0,1)$, so that the generating curves $\Gamma(t)$ are open, and the
surfaces $\mathcal{S}(t)$ have genus 0. The corresponding scheme for the case
$I=\RZ$ can then easily be obtained by ignoring the special treatment of
the boundary nodes.

In order to define our finite difference schemes, let us introduce the set of 
grid 
points $\mathcal G^h:= \{\rho_0,\ldots, \rho_J \} \subset I$, where 
$\rho_j = jh$, $j=0,\ldots, J$, and $h = \frac 1J$ for $J\geq2$. 
We also make use of the same time discretization as before,
with $t_m=m \Delta t$.
For a grid function $v^m:\mathcal G^h \to \bR^2$ we write 
$v^m_j:= v^m(\rho_j)$, $j=1,\ldots,J$. 
We associate with $v^m$ the difference quotients
\begin{alignat*}{2}
\delta^- v^m_j &:= \frac{v^m_j - v^m_{j-1}}{h}, \quad &&j=1,\ldots,J; \\
\delta^+ v^m_j &:= \frac{v^m_{j+1} - v^m_j}{h}, \quad &&
j=0,\ldots,J-1; \\
\delta^1 v^m_j &:= \frac{v^m_{j+1} - v^m_{j-1}}{2h},\quad &&
j=1,\ldots,J-1. 
\end{alignat*}

Let $x^m: \mathcal G^h \to \bR^2$ be a grid function that approximates
$x(\cdot,t_m)$. Then on $I_j=[\rho_{j-1},\rho_j]$, 
the associated discrete length element $q^m_j$ and the discrete tangent 
$\tau^m_j$ are given by 
\begin{equation*} 
q^m_j = | \delta^- x^m_j |, \quad \tau^m_{j} = \frac{1}{q^m_j} \delta^- x^m_j, 
\quad j=1,\ldots,J.
\end{equation*}
It will be convenient to also introduce the averaged vertex tangent 
$\theta^m_j$ via
\begin{equation} \label{eq:app:fdthetah}
\theta^m_j = 
\frac{\tau^m_{j} + \tau^m_{j+1}}{|\tau^h_{j} + \tau^m_{j+1}|}, \quad \mbox{ provided that } \tau^m_j + \tau^m_{j+1} \neq 0, \quad j=1,\ldots,J-1.
\end{equation}

\begin{remark} \label{rem:hmcf}
We observe that \eqref{eq:app:hcsfa} is very close to
the formulation for hyperbolic curvature flow for closed curves
used in \cite{hmcf}: 
\begin{equation} \label{eq:hmcf}
x_{tt} = \kappa \nu - (x_t \cdot \tau_t) \tau
\quad \mbox{ in } I \times (0,T].
\end{equation}
There the authors presented the following fully discrete finite difference 
approximation for \eqref{eq:hmcf}.
Given suitable initial data $x^0, x^{-1}: \mathcal G^h \to \bR^2$,
for $m=0,\ldots,M-1$ we find $x^{m+1}: \mathcal G^h \to \bR^2$ such 
that, for $j = 1,\ldots,J$, 
\begin{align} \label{eq:app:fdfdimpl}
& \tfrac12 (q^m_j + q^m_{j+1}) h
\frac{ x^{m+1}_j - 2  x^m_j +  x^{m-1}_j}{(\Delta t)^2} 
\nonumber \\ & \quad 
= 
\frac{1}{2q^m_{j+1}} \delta^- x^{m+1}_{j+1} - \frac{1}{2q^m_j} \delta^- x^{m+1}_j
+\frac{1}{2q^m_{j+1}} \delta^- x^{m-1}_{j+1} - \frac{1}{2q^m_j} \delta^- x^{m-1}_j
\nonumber \\ & \qquad 
- \tfrac12 (q^m_j + q^m_{j+1}) h
(\frac{x^m_j -  x^{m-1}_j}{\Delta t} \cdot \frac{\theta^m_j - \theta^{m-1}_j}
{\Delta t}) \theta^m_j. 
\end{align}
We also recall for the reader that the choice of $\theta^m_j$ in
\eqref{eq:app:fdfdimpl}, as defined in \eqref{eq:app:fdthetah}, 
is motivated by the desire to prove a semidiscrete 
analogue of the normal parameterization property, $x_t \cdot \tau =0$.
In fact, in Lemma~3.2 in \cite{hmcf} it was shown that the continuous-in-time
semidiscrete variant of \eqref{eq:app:fdfdimpl} 
satisfies such a discrete analogue.
\end{remark}

Following \eqref{eq:app:fdfdimpl}, comparing \eqref{eq:hmcf} with 
\eqref{eq:app:hcsfa}, and noting that
\begin{equation} \label{eq:mcax}
\kappa \nu - \frac{\nu \cdot e_1}{x \cdot e_1} \nu =
\frac{1}{|x_\rho|} \bigl(\frac{x_\rho}{|x_\rho|}\bigr)_\rho - \frac{1}{|x_\rho|} 
\frac{x_\rho \cdot e_2}{x \cdot e_1} \tau^\perp 
\quad \mbox{ in } I \times (0,T],
\end{equation}
we propose the following fully discrete 
approximation for the flow \eqref{eq:app:hcsf}. 
Given suitable initial data $x^0, x^{-1}: \mathcal G^h \to \bR^2$,
for $m=0,\ldots,M-1$ we find $x^{m+1}: \mathcal G^h \to \bR^2$ such 
that 
\begin{subequations} \label{eq:app:fd}
\begin{align} \label{eq:app:fda}
& \tfrac12 (q^m_j + q^m_{j+1}) h
\frac{ x^{m+1}_j - 2  x^m_j +  x^{m-1}_j}{(\Delta t)^2} \nonumber \\
&
=  \hat g^m_j \left[\frac{1}{2q^m_{j+1}} \delta^- x^{m+1}_{j+1} - \frac{1}{2q^m_j} \delta^- x^{m+1}_j 
+\frac{1}{2q^m_{j+1}} \delta^- x^{m-1}_{j+1} - \frac{1}{2q^m_j} \delta^- x^{m-1}_j
- \frac{\delta^1  x^{m}_j \cdot  e_2 }{x^m_j \cdot e_1}
(\theta^m_j)^\perp \right]
\nonumber \\ & \qquad 
- \tfrac12 (q^m_j + q^m_{j+1}) h
(\frac{x^m_j -  x^{m-1}_j}{\Delta t} \cdot \frac{\theta^m_j - \theta^{m-1}_j}
{\Delta t}) \theta^m_j, \quad j = 1,\ldots,J-1, \\
& x^{m+1}_0 \cdot e_1 = 0; \quad 
\delta^+ x^{m+1}_0 \cdot e_2 = \tfrac14 h  \frac{( q^m_1 )^2}{\hat g^m_0 }
 \frac{x^{m+1}_0 - 2 x^m_0 + x^{m-1}_0}{(\Delta t)^2} 
 \cdot e_2 , 
\label{eq:app:fdb}  \\
& x^{m+1}_J \cdot e_1 = 0; \quad 
\delta^- x^{m+1}_J \cdot e_2 = - \tfrac14 h \frac{( q^m_J )^2}{ \hat g^m_J }
 \frac{x^{m+1}_J - 2 x^m_J + x^{m-1}_J}{(\Delta t)^2} 
 \cdot e_2,
\label{eq:app:fdc}
\end{align}
\end{subequations}
where $\hat g^m_j = g ( |\frac{x^m_j -  x^{m-1}_j}{\Delta t}|^2 )$. 
In order to motivate the second conditions in \eqref{eq:app:fdb} and \eqref{eq:app:fdc}, we first calculate with the help of a Taylor expansion and \eqref{eq:app:hcsfb} that
\begin{equation} \label{eq:cont:fdb}
\frac{x(h,t_{m+1})-x(0,t_{m+1})}{h} \cdot e_2 \approx x_\rho(0,t_{m+1}) \cdot e_2 + \tfrac12 h x_{\rho \rho}(0,t_m) \cdot e_2 = \tfrac12 h x_{\rho \rho}(0,t_m) \cdot e_2.
\end{equation}
Moreover, we deduce from  \eqref{eq:mcax}, \cite[(A3)]{degenD} and the identity $x_\rho(0,t_m)= | x_\rho(0,t_m) | e_1$ that
\begin{align*} 
 \lim_{\rho \searrow 0} \bigl( \kappa \nu - \frac{\nu \cdot e_1}{x \cdot e_1} \nu \bigr)(\rho,t_m)  & = 
 \left[ \frac{x_{\rho \rho}}{| x_\rho|^2} - \frac{x_{\rho \rho} \cdot x_\rho}{| x_\rho |^4} x_\rho \right](0,t_m) -  \lim_{\rho \searrow 0} \left[\frac{1}{|x_\rho|^2} 
\frac{x_\rho \cdot e_2}{x \cdot e_1} x_\rho^\perp \right](\rho,t_m) \\
& = 2 \frac{x_{\rho \rho}(0,t_m) \cdot e_2}{| x_\rho(0,t_m) |^2} e_2.
\end{align*}
In addition, we see from \eqref{eq:app:hcsfb} that
$x_t(0,t_m) \cdot e_1 = 0$ and $\tau_t(0,t_m) \cdot e_2=0$. 
Hence $(x_t \cdot \tau_t)(0,t_m)=0$, so that \eqref{eq:app:hcsfa} yields
\begin{displaymath}
x_{tt}(0,t_m) \cdot e_2 = g(|x_t(0,t_m)|^2) e_2 \cdot \lim_{\rho \searrow 0} \bigl( \kappa \nu - \frac{\nu \cdot e_1}{x \cdot e_1} \nu \bigr)(\rho,t_m) = 2 g(|x_t(0,t_m)|^2) \frac{x_{\rho \rho}(0,t_m) \cdot e_2}{| x_\rho(0,t_m) |^2} ,
\end{displaymath}
or equivalently
\begin{displaymath}
\tfrac12 h x_{\rho \rho}(0,t_m) \cdot e_2 = \tfrac14 h \frac{| x_\rho(0,t_m)|^2}{ g(| x_t(0,t_m)|^2)} x_{tt}(0,t_m) \cdot e_2.
\end{displaymath}
Inserting this relation into \eqref{eq:cont:fdb} and replacing derivatives by suitable difference quotients, we are led to \eqref{eq:app:fdb}. A similar argument applies to \eqref{eq:app:fdc}.

\begin{lemma} \label{lem:axiex}
Let $g \in C^0(\bRgeq)$ with $g(s) \geq0$. 
Suppose that $q^m_j>0$ for $j=1,\ldots,J$. Then there exists a unique solution $x^{m+1}_h : \mathcal{G}^h \to \bR^2$
to \eqref{eq:app:fd}. 
\end{lemma}
\begin{proof} 
Since \eqref{eq:app:fd} is a linear scheme with the same number of unknowns as
equations, it is sufficient to prove that only the trivial solution solves
the homogeneous system
\begin{subequations} \label{eq:proof:fd}
\begin{align} \label{eq:proof:fda}
& \tfrac12 (q^m_j + q^m_{j+1}) h^2
\frac{ X_j }{(\Delta t)^2} 
= 
\frac{\hat g^m_j}{2q^m_{j+1}} (X_{j+1} - X_j) 
- \frac{\hat g^m_j}{2q^m_j} (X_j - X_{j-1})
, \quad j = 1,\ldots,J-1,\\
& X_0 \cdot e_1 = 0; \quad 
(X_1 - X_0) \cdot e_2 = \tfrac14 h^2 \frac{( q^m_1 )^2}{\hat g^m_0 }
\frac{X_0}{(\Delta t)^2} \cdot e_2 , 
\label{eq:proof:fdb}  \\
& X_J \cdot e_1 = 0; \quad 
(X_J -X_{J-1}) \cdot e_2 = - \tfrac14 h^2 \frac{ ( q^m_J )^2}{ \hat g^m_J }
\frac{X_J }{(\Delta t)^2} \cdot e_2.
\label{eq:proof:fdc}
\end{align}
\end{subequations}
Clearly, \eqref{eq:proof:fd} decouples into two tridiagonal linear systems
for $X \cdot e_1$ and $X \cdot e_2$, which are both strictly diagonally
dominant. Hence it follows that $X=0$.
\end{proof}

In order to derive suitable initial data for \eqref{eq:app:fd}, 
and on recalling that for simplicity we assume $V_0$ to be constant,
we again use a formal Taylor expansion together with \eqref{eq:app:hcsf}
to calculate in $(0,1)$
\begin{align} \label{eq:app:x0Taylor}
x(\cdot, -\Delta t) & = x(\cdot,0) - \Delta t x_t (\cdot,0)
+ \tfrac12 (\Delta t)^2 x_{tt}(\cdot,0) + \mathcal{O}((\Delta t)^3)
\nonumber \\ &
= x_0 - \Delta t V_0 \nu(\cdot,0) 
+ \tfrac12 (\Delta t)^2 \left[ g(V_0^2) \Bigl( \frac{1}{| x_\rho |}
\bigl(\frac{x_\rho}{| x_\rho |} \bigr)_\rho 
- \frac{\nu \cdot e_1}{x \cdot e_1} \nu  \Bigr)
- ( V_0 \nu \cdot \tau_t) \tau \right](\cdot,0)
\nonumber \\ & \qquad
+ \mathcal{O}((\Delta t)^3)
\nonumber \\ &
= x_0- \Delta t V_0 \nu(\cdot,0)
+ \tfrac12 (\Delta t)^2 g(V_0^2) \left[ \frac{1}{| x_\rho |}
\bigl(\frac{x_\rho}{| x_\rho |} \bigr)_\rho 
- \frac{\nu \cdot e_1}{x \cdot e_1} \nu 
\right](\cdot,0)
+ \mathcal{O}((\Delta t)^3),
\end{align}
where we used in the last step that at time $t=0$
\begin{displaymath}
\tau_t \cdot \nu = \frac{x_{t, \rho}}{| x_\rho |} \cdot \nu = \frac{1}{| x_\rho|} (V_0 \nu)_\rho \cdot \nu = \frac{1}{| x_\rho |} V_0 \nu_\rho \cdot \nu=0.
\end{displaymath}
In addition we have from \cite[(2.17)]{axisd} that
\begin{displaymath}
\lim_{\rho\to\rho_0} \frac{\nu \cdot e_1}{x \cdot e_1} \nu = 
-\kappa(\rho_0)\nu(\rho_0) \qquad \forall\ \rho_0 \in \{0,1\},
\end{displaymath}
so that 
\begin{equation} \label{eq:app:kappa0}
x(\rho_0,-\Delta t) = x_0(\rho_0) - \Delta t V_0 \nu(\rho_0,0) + 
\tfrac12 (\Delta t)^2 g(V_0^2) [2\kappa \nu](\rho_0,0) 
+ \mathcal{O}((\Delta t)^3)
\quad \forall\ \rho_0 \in \{0,1\}.
\end{equation}
Inspired by \eqref{eq:app:x0Taylor} and \eqref{eq:app:kappa0},
we choose as initial data
\begin{alignat*}{2} 
x^0_j & = x_0(\rho_j)\ \qquad && j = 0,\ldots,J,\\
x^{-1}_j & = x^0_j - \Delta t V_0 \theta^{0,\perp}_j 
+ \tfrac12 (\Delta t)^2 g(V_0^2) 
y^0_j ,
\quad && j=1,\ldots,J-1, \\ 
x^{-1}_j & = x^0_j + \Delta t V_0 \sign(\tau^0_j \cdot e_1)
e_2 + \tfrac12 (\Delta t)^2 g(V_0^2) y^0_j , \quad && j=0,J,
\end{alignat*}
where we have defined $\tau^0_0 := \tau^0_1$.
Here $y^m:\mathcal{G}^h \to \bR^2$ is a grid function that represents a
discrete surface mean curvature vector defined as follows:
\begin{alignat}{2}
y^m_j &= \frac{1}{q^m_{j+1}} \delta^- x^m_{j+1} - \frac{1}{q^m_j} \delta^- x^m_j
- \frac{\delta^1  x^m_j \cdot  e_2 }{x^m_j \cdot e_1}
(\theta^m_j)^\perp, \quad && j = 1,\ldots,J-1, \nonumber \\
y^m_j &= \frac{2}{q^m_{j+1}} \delta^- x^m_{j+1} - \frac{2}{q^m_j} \delta^- x^m_j, \quad && j = 0,J,
\label{eq:kappa0}
\end{alignat}
where here we have introduced the ``ghost points'' $x^m_{-1}$ and $x^m_{J+1}$
that are obtained by a reflection on the $x_2$-axis of $x^m_1$ and $x^m_{J-1}$,
respectively.

\setcounter{equation}{0}
\section{Numerical results} \label{sec:nr}

We implemented all our numerical schemes within the
finite element toolbox Alberta, \cite{Alberta}, using
the sparse factorization package UMFPACK, see \cite{Davis04},
for the solution of the linear systems of equations arising at each time 
level.

A particularly fascinating aspect of hyperbolic mean curvature flow is the
possible development of singularities in finite time from smooth initial data.
As observed in \cite{hmcf} for the corresponding flow for curves in the plane,
some of these singularities may be characterized by a blow-up in curvature due
to the surface developing kinks or sharp corners. In order to monitor the onset
of such singularities, we will often be interested in the quantity
\[
\mathcal K^m_\infty = \max_{j=0,\ldots,J} |y^m_j|,
\]
for solutions of the scheme \eqref{eq:app:fd},
were $y^m:\mathcal{G}^h \to \bR^2$ is defined by \eqref{eq:kappa0}. In
addition, we define the surface area of the surfaces generated by the solutions
to \eqref{eq:app:fd} via
\[
A^m = 2\pi h \sum_{j=1}^J x^m_j \cdot e_1 |\delta^- x^m_j|.
\]

We remark that in all our computations for \eqref{eq:Gurtinbeta0}
and \eqref{eq:LeFloch}, respectively, natural discrete analogues of
$E(t)= \tfrac12 \int_{\mathcal S(t)} (\mathcal{V}_{\mathcal S}^2 + 2) \dH2$
and $\widetilde E(t)= \int_{\mathcal S(t)} e^{\frac12\mathcal{V}_{\mathcal S}^2} \dH2$ 
were well preserved, recall \eqref{eq:Gurtinbeta0} and \eqref{eq:ELeFloch}. 
In particular, in all our computations the relative
differences between these discrete analogues at time $t=0$ and at time $t=T$ 
are always less than $1.5\%$, even though the computations often come very
close to a singularity.

\subsection{The hyperbolic mean curvature flow \eqref{eq:Gurtinbeta0}}

\subsubsection{Radial evolutions for a sphere}

Our first set of numerical experiments is for the evolution of a sphere.
To validate our schemes \eqref{eq:fea_b} and \eqref{eq:app:fd} with $g(s)=1$, 
we begin with a
convergence experiment, generalizing the solution from \cite[\S5.2]{hmcf} 
to higher space dimensions. It can be shown that a family of spheres in $\mathbb R^3$ with radius $r(t)$ and outer
normal $\nu(t)$ is a solution to \eqref{eq:Gurtinbeta0} with 
$\mathcal V_{\mathcal{S} \mid_{t=0}} = V_0 \in \bR$ if
\begin{equation*} 
\ddot{r}(t) = - \frac{2}{r(t)} \quad \text{in } (0,T],\quad
r(0) = r_0>0,\ \dot{r}(0) = V_0.
\end{equation*}  
Upon integration we obtain that
\[
\tfrac12 (\dot{r}(t))^2 = 2 (\ln r_0 - \ln r(t)) + \tfrac12 V_0^2
= 2 \ln \frac{r_0}{r(t)} + \tfrac12 V_0^2.
\]
Hence
\begin{equation} \label{eq:rdot}
\dot{r}(t) = \pm \sqrt{4 \ln \frac{r_0}{r(t)} + V_0^2 },
\end{equation}
which means that if $V_0 > 0$, then $r(t)$ will at first increase until it hits
a maximum, where $4 \ln \frac{r_0}{r(t)} + V_0^2 = 0$, after which it
will decrease and shrink to a point in finite time. On the other hand,
if $V_0 \leq 0$ then the sphere will monotonically shrink to a point.

For the special case $V_0 = 0$, and on recalling the Gauss error function
\begin{equation*} 
\erf(z) = \tfrac2{\sqrt\pi} \int_0^z e^{-u^2} \;{\rm d}u,\quad
\erf'(z) = \tfrac2{\sqrt\pi} e^{-z^2},
\end{equation*}
we find that $r(t)$ is the solution of
\begin{equation*} 
t - \sqrt{\tfrac\pi{4}} r_0 \erf\left(\sqrt{\ln\frac{r_0}{r(t)}}\right) = 0,
\end{equation*}
which means that
\begin{equation} \label{eq:rsol}
r(t) = r_0 \exp(-[\erf^{-1}(\sqrt{\tfrac{4}\pi}\frac{t}{r_0})]^2).
\end{equation}

We now use the exact solution \eqref{eq:rsol} in order to validate the scheme \eqref{eq:fea_b}.
Using a sequence of triangulations of the unit sphere yields the results
reported in Table~\ref{tab:3d}, where we have defined the error
\begin{equation} \label{eq:errorSS}
\errorSS = \max_{m=1,\ldots,M} \max_{k=1,\ldots,K} | |p^m_k| - r(t_m)|,
\end{equation}
and where we employed the time step size choice $\Delta t = 0.25 h_0$.
We note that the observed order of convergence
approaches the optimal rate $\mathcal{O}(h^2_0)$.

\begin{table}
\center
\begin{tabular}{|r|c|c|c|}
\hline 
$J$   & $h_0$      & $\errorSS$ & EOC \\ \hline
1536  & 2.0854e-01 & 5.0520e-03 & --- \\
6144  & 1.0472e-01 & 1.5715e-03 & 1.70 \\
24576 & 5.2416e-02 & 4.5239e-04 & 1.80 \\ 
98304 & 2.6215e-02 & 1.2384e-04 & 1.87\\
393216& 1.3108e-02 & 3.3629e-05 & 1.88\\
\hline
\end{tabular}
\caption{Errors for \eqref{eq:fea_b} for the convergence test for 
\eqref{eq:rsol} with $r_0 = 1$ over the time interval $[0,0.25]$.
We also display the experimental orders of convergence (EOC).
}
\label{tab:3d}
\end{table}%

Next, we perform the same convergence experiment also for the axisymmetric
scheme \eqref{eq:app:fd} with $g(s)=1$. To this end, we choose as initial data
$x_0(\rho) = (\cos (\frac32\pi + \pi \rho), \sin (\frac32\pi + \pi \rho) )^T$. It is easily verified that $x(\rho,t)=r(t) x_0(\rho)$ is an exact solution
of \eqref{eq:app:hcsf} with $g(s) = 1$. 
The observed errors are shown in Table~\ref{tab:axi}, where we have
defined 
\begin{equation} \label{eq:errorXX}
\errorXX = \max_{m=1,\ldots,M} \max_{j=0,\ldots,J} | x(\rho_j,t_m) - x^m_j |,
\end{equation} 
and where we have chosen $\Delta t =h ~= \frac 1J$.
Also for the scheme \eqref{eq:app:fd} we obtain an optimal convergence rate.
\begin{table}
\center
\begin{tabular}{|r|c|c|}
\hline
$J$ & $\errorXX$ & EOC \\ \hline
32  & 6.3402e-04 & ---  \\
64  & 1.3346e-04 & 2.25 \\
128 & 2.9262e-05 & 2.19 \\ 
256 & 6.9967e-06 & 2.06 \\ 
512 & 1.7053e-06 & 2.04 \\ 
\hline
\end{tabular}
\caption{Errors for \eqref{eq:app:fd} with $g(s) = 1$ 
for the convergence test for 
\eqref{eq:rsol} with $r_0=1$ over the time interval $[0,0.5]$. 
We also display the experimental orders of convergence (EOC).
}
\label{tab:axi}
\end{table}%

\begin{figure}
\center
\includegraphics[angle=-90,width=0.3\textwidth]{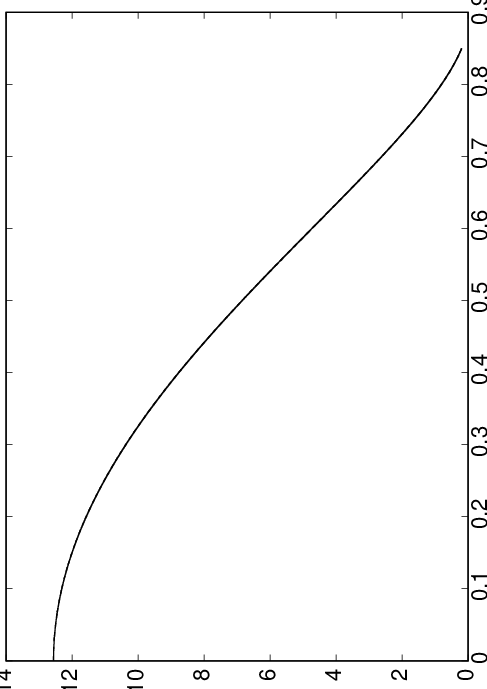}
\includegraphics[angle=-90,width=0.3\textwidth]{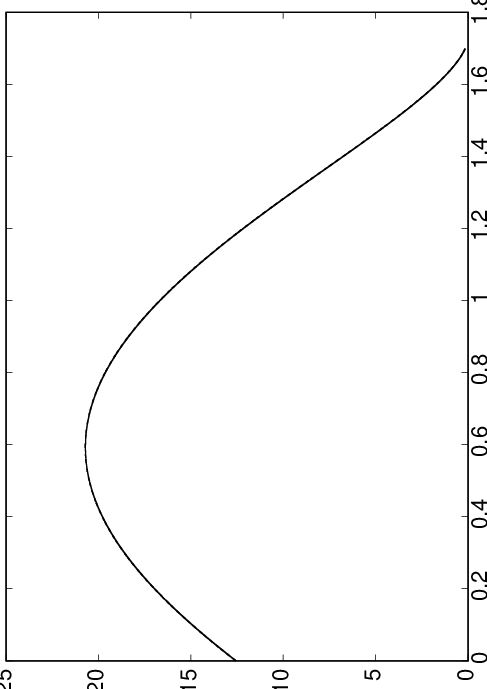}
\includegraphics[angle=-90,width=0.3\textwidth]{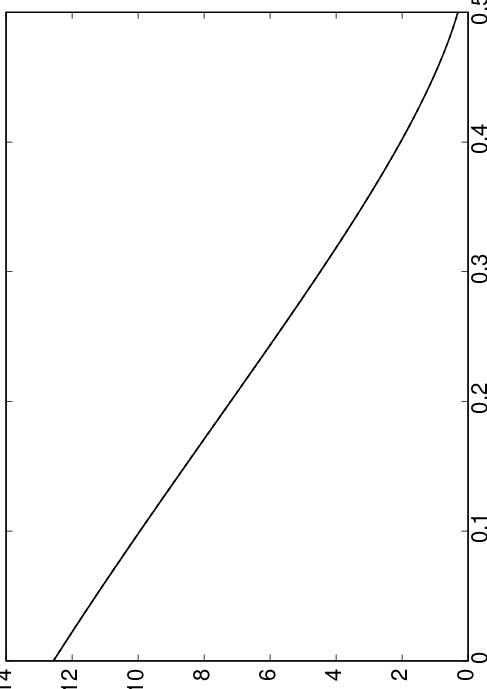}
\\
\includegraphics[angle=-90,width=0.3\textwidth]{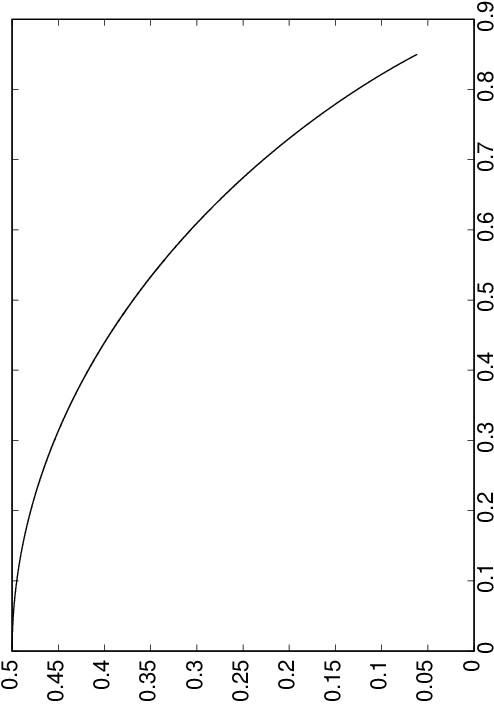}
\includegraphics[angle=-90,width=0.3\textwidth]{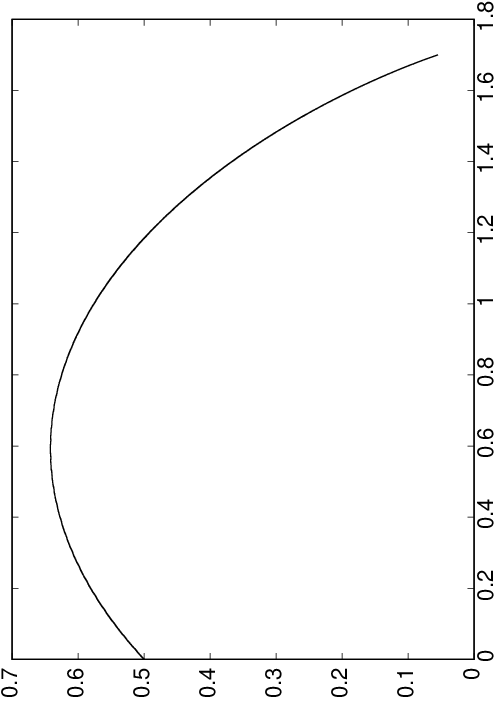}
\includegraphics[angle=-90,width=0.3\textwidth]{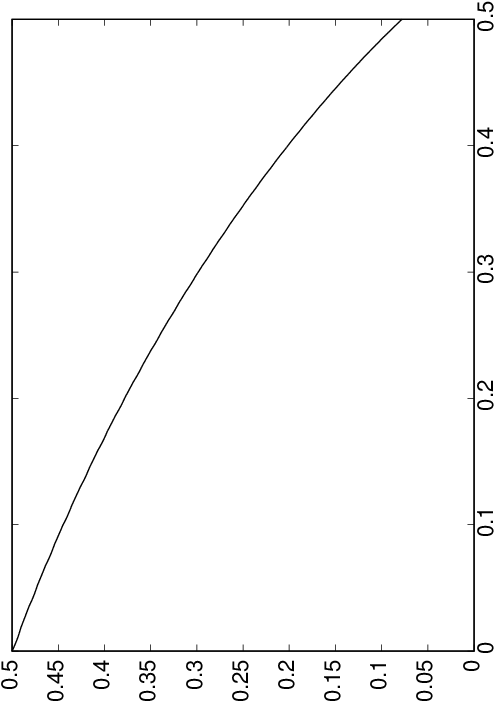}
\caption{The flow \eqref{eq:Gurtinbeta0} starting from a unit sphere.
Above we show the evolutions of $|\mathcal S^m_h|$ (solid) and $A^m$ (dashed) 
over time for $V_0=0$, $V_0=1$ and $V_0=-1$.
The final times for these computations are $T=0.85$, $T=1.7$ and $T=0.5$,
respectively.
Below we show the corresponding evolutions of $1/{\mathcal K}^m_\infty$ over 
time.
}
\label{fig:axispheres}
\end{figure}

As discussed above, for an initial sphere, depending on the sign of 
$V_0 \in \bR$, the family of
spheres either expands at first and then shrinks, or shrinks immediately. We
visualize these different behaviours in Figure~\ref{fig:axispheres}. 
Here we plot
$|\mathcal S^m_h|$ and $A^m$ in the same figure, to confirm perfect agreement
between our two approximations \eqref{eq:fea_b} and \eqref{eq:app:fd}
with $g(s)=1$. 
In fact here and in what follows the plots of the two quantities lie on top
of each other and so are indistinguishable to the human eye.
In each case we observe a smooth solution until the spheres shrink to a point, 
meaning that $|\mathcal S^m_h|$ and $1/{\mathcal K}^m_\infty$ approach zero at the same time.

\subsubsection{Axisymmetric evolutions for nonspherical surfaces}
\label{sec:nraxi}

For all the numerical simulations in this subsection we 
choose the time discretization parameter $\Delta t = 10^{-4}$. 

\begin{figure}
\center
\includegraphics[angle=-90,width=0.15\textwidth]{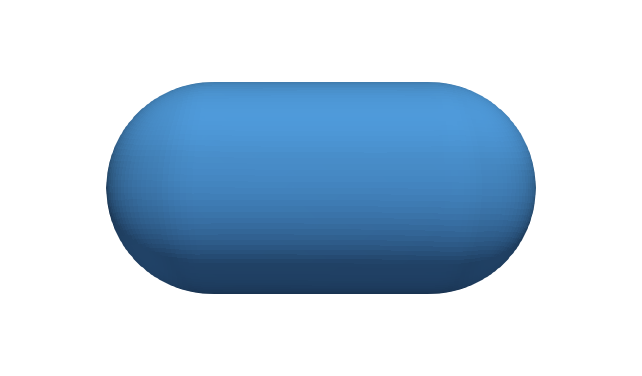} 
\includegraphics[angle=-90,width=0.15\textwidth]{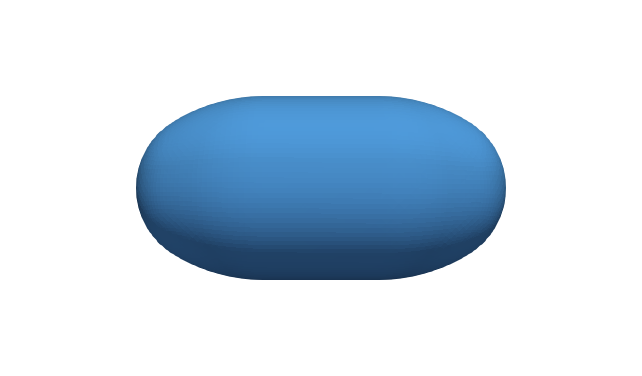} 
\includegraphics[angle=-90,width=0.15\textwidth]{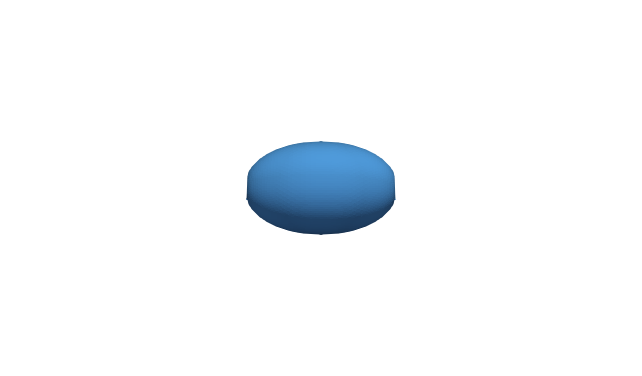} 
\includegraphics[angle=-90,width=0.15\textwidth]{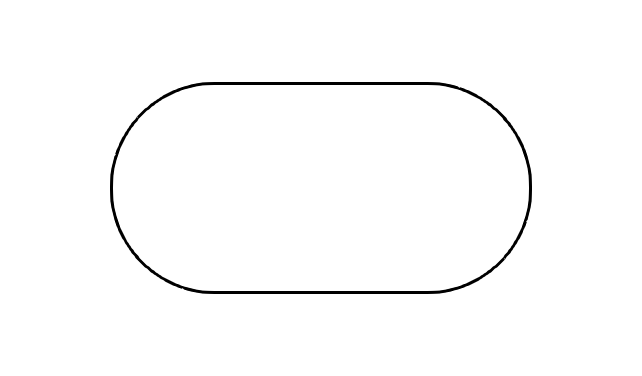} 
\includegraphics[angle=-90,width=0.15\textwidth]{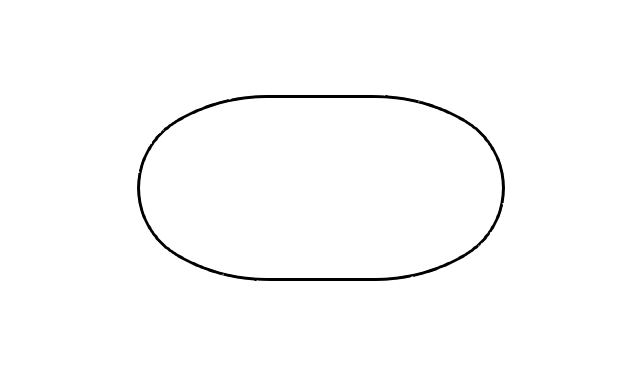} 
\includegraphics[angle=-90,width=0.15\textwidth]{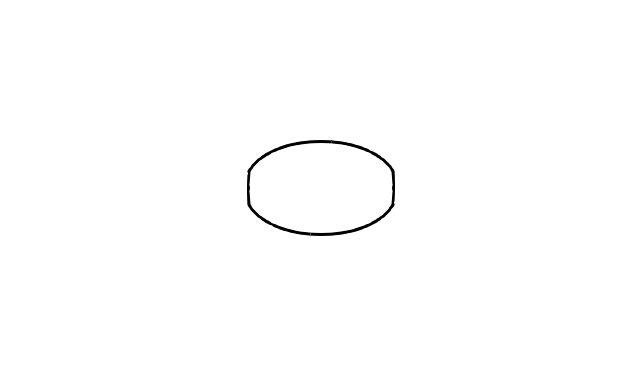} 
\\
\center
\includegraphics[angle=-90,width=0.12\textwidth]{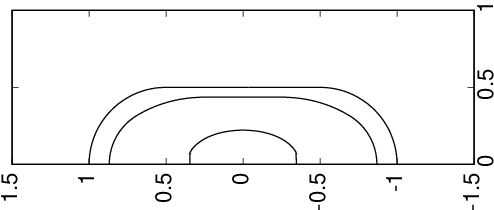}
\includegraphics[angle=-90,width=0.4\textwidth]{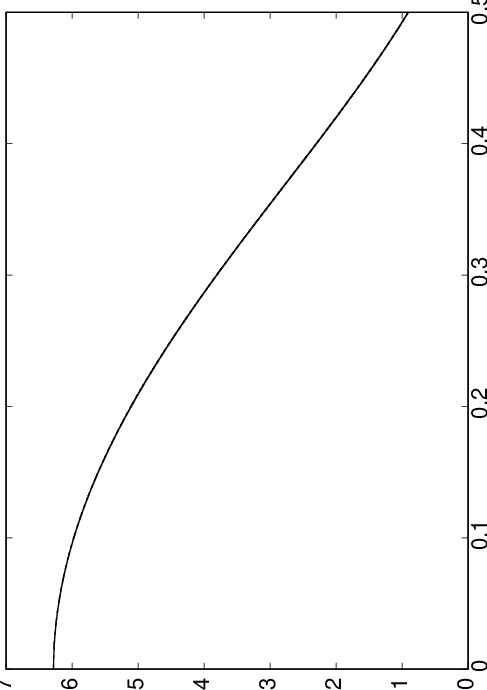}
\includegraphics[angle=-90,width=0.4\textwidth]{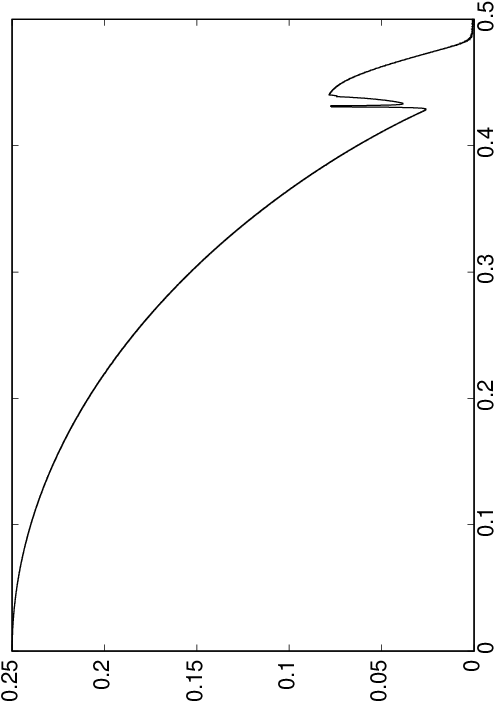}
\caption{The flow \eqref{eq:Gurtinbeta0}, with $V_0=0$, starting from a 
$2:1:1$ cigar.
Above we visualize $\mathcal S^m_h$ from \eqref{eq:fea_b} 
at times $t=0,0.25,T=0.5$ in two ways: a frontal view and a cut through the
$x_1$-$x_2$ plane.
Below we show $x^m$ from \eqref{eq:app:fd} with $g(s)=1$
at times $t=0,0.25,T$, as well as
a plot with $|\mathcal S^m_h|$ (solid) and $A^m$ (dashed), and a plot of
$1/\mathcal{K}^m_\infty$ (right) over time.
}
\label{fig:cigar211V0}
\end{figure}

We begin with an experiment for an initial $2:1:1$ cigar,
when the initial velocity is chosen as $V_0=0$,
see Figure~\ref{fig:cigar211V0}. 
Again we observe an excellent agreement between the solutions for our schemes
\eqref{eq:fea_b} and \eqref{eq:app:fd} with $g(s)=1$.
At around  $t=0.5$ there seems to be the onset of a singularity.
Notice that the first kink in the curvature plot 
at around time $t=0.43$ indicates a temporary increase in the
magnitude of the mean 
curvature at the axis of rotation. However, this does not appear to lead to a
singularity. Instead, the shorter part of the surface eventually becomes flat,
leading to a blow-up of curvature where the flat part meets the curved
part of the surface. To qualitatively confirm these results, we have repeated
this simulation also with finer discretization parameters, and the described
behaviour is consistent across the different computations.

\begin{figure}
\center
\includegraphics[angle=-90,width=0.15\textwidth]{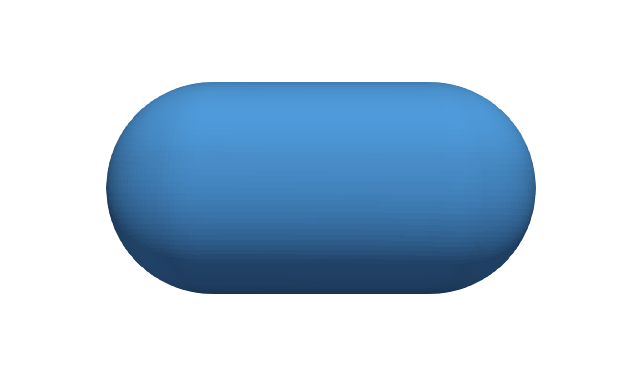} 
\includegraphics[angle=-90,width=0.15\textwidth]{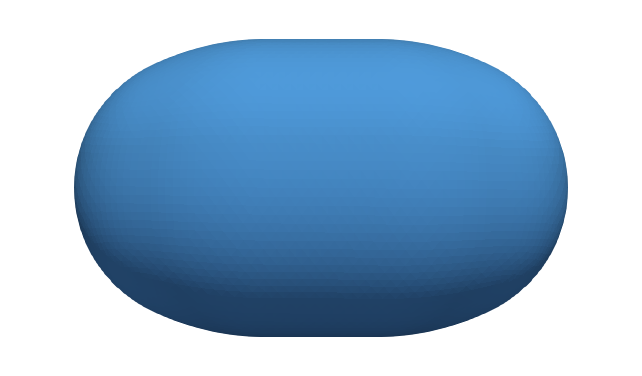} 
\includegraphics[angle=-90,width=0.15\textwidth]{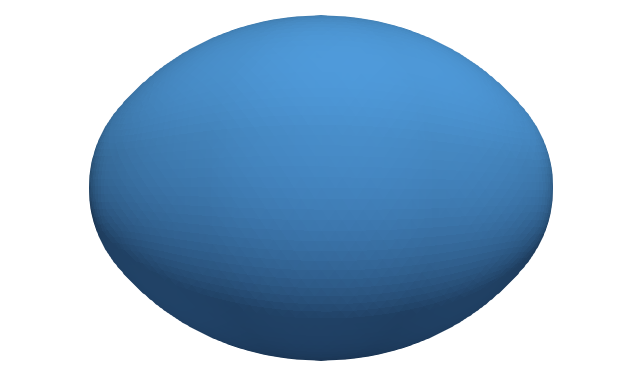} 
\includegraphics[angle=-90,width=0.15\textwidth]{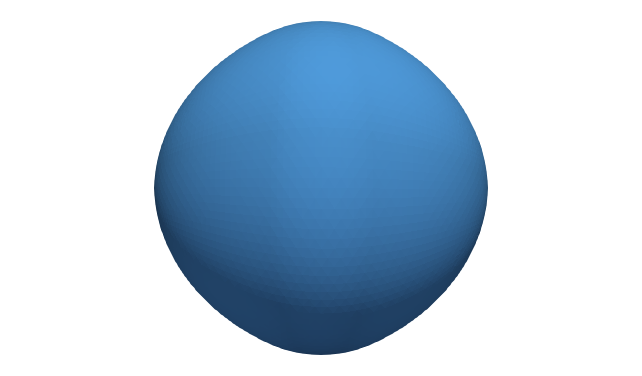} 
\includegraphics[angle=-90,width=0.15\textwidth]{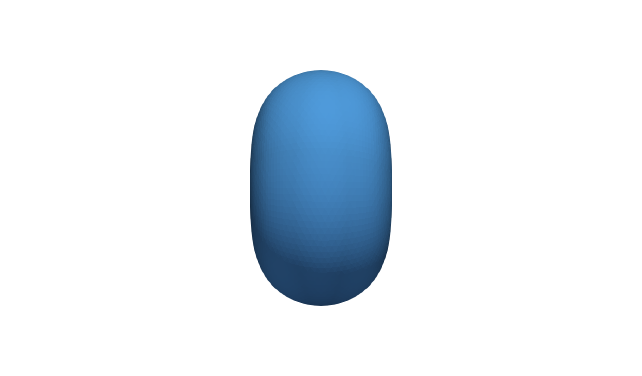} 
\includegraphics[angle=-90,width=0.15\textwidth]{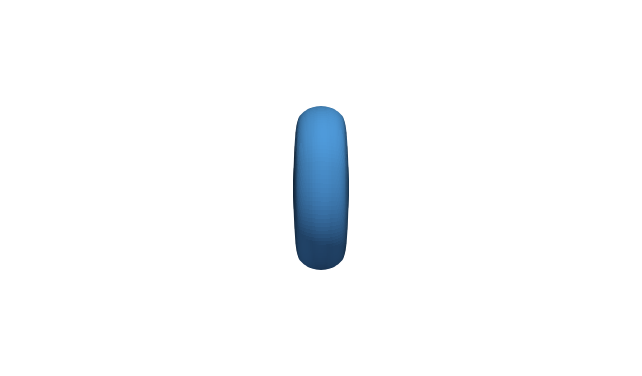} 
\includegraphics[angle=-90,width=0.15\textwidth]{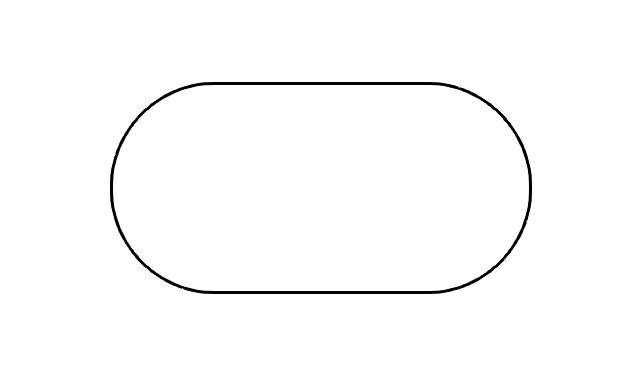} 
\includegraphics[angle=-90,width=0.15\textwidth]{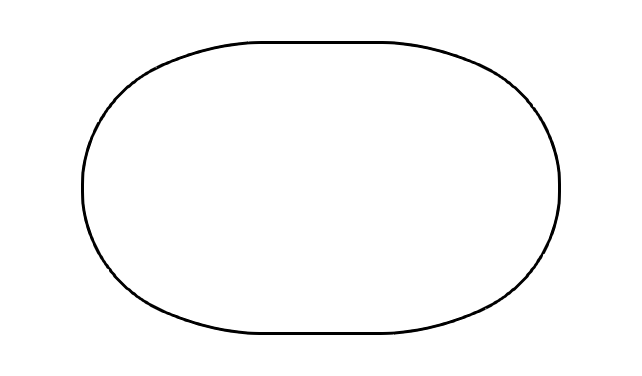} 
\includegraphics[angle=-90,width=0.15\textwidth]{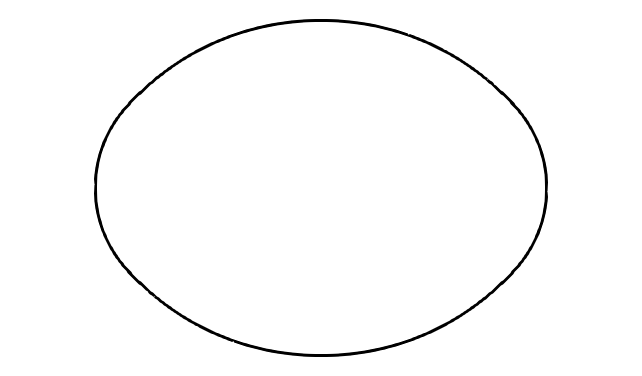} 
\includegraphics[angle=-90,width=0.15\textwidth]{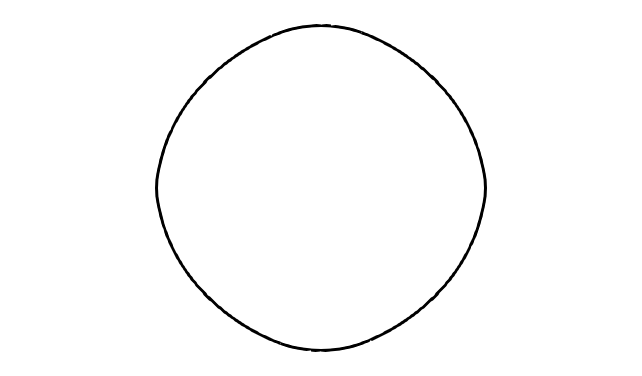} 
\includegraphics[angle=-90,width=0.15\textwidth]{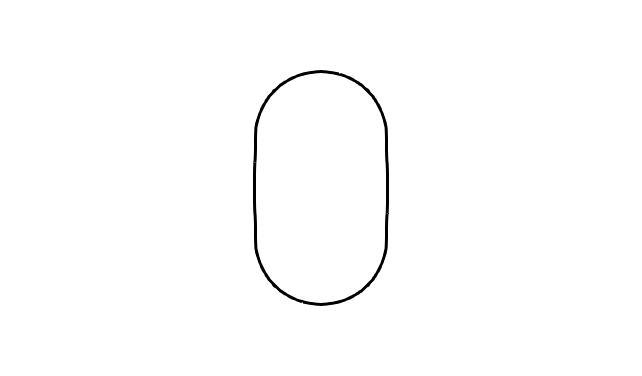} 
\includegraphics[angle=-90,width=0.15\textwidth]{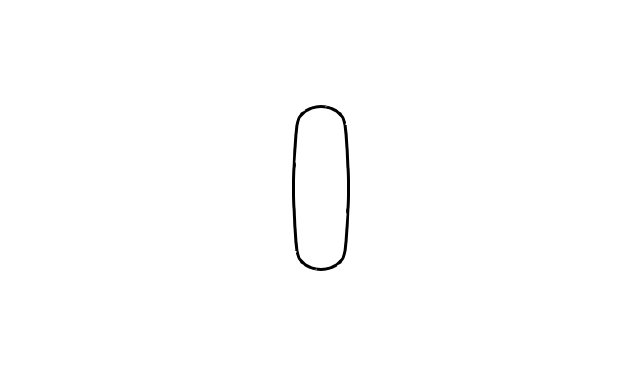} 
\\
\center
\includegraphics[angle=-90,width=0.12\textwidth]{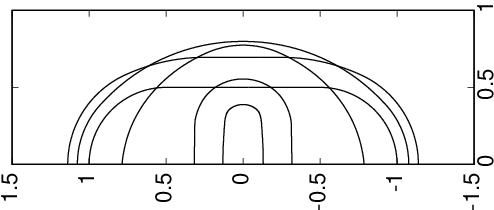}
\includegraphics[angle=-90,width=0.4\textwidth]{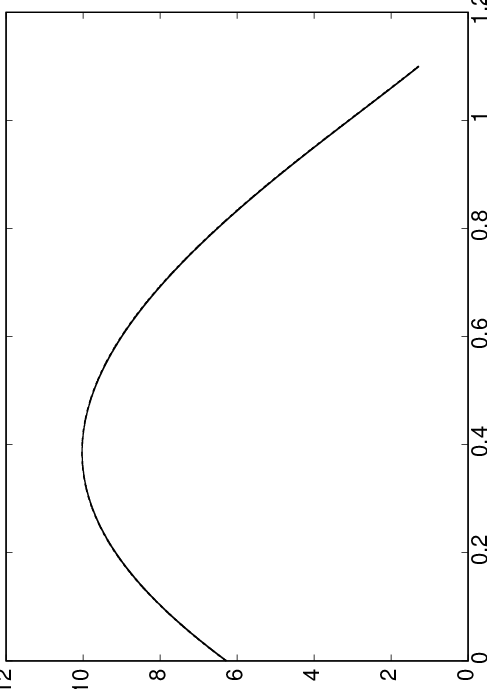}
\includegraphics[angle=-90,width=0.4\textwidth]{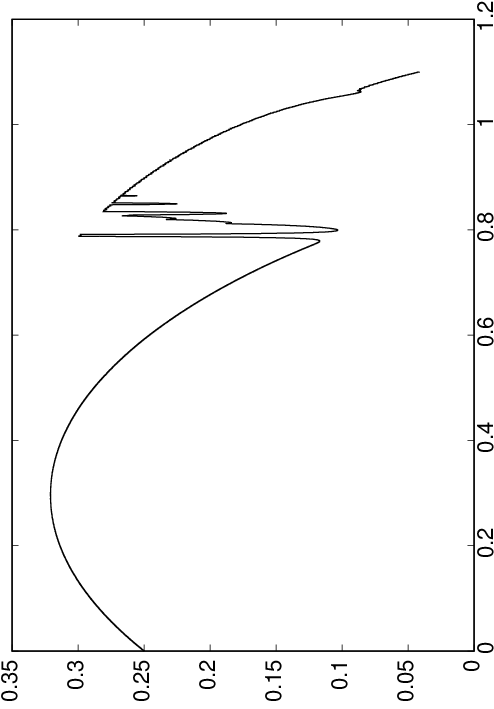}
\caption{The flow \eqref{eq:Gurtinbeta0}, with $V_0=1$, starting from a 
$2:1:1$ cigar.
Above we visualize $\mathcal S^m_h$ from \eqref{eq:fea_b} 
at times $t=0,0.25,0.5,0.75,1,T=1.1$ in two ways: 
a frontal view and a cut through the $x_1$-$x_2$ plane.
Below we show $x^m$ from \eqref{eq:app:fd} with $g(s)=1$
at times $t=0,0.25,0.5,0.75,1,T$, as well as
a plot with $|\mathcal S^m_h|$ (solid) and $A^m$ (dashed), and a plot of
$1/\mathcal{K}^m_\infty$ (right) over time.
}
\label{fig:cigar211V1}
\end{figure}

When we change the initial velocity to $V_0=1$, the evolution is dramatically
different, see Figure~\ref{fig:cigar211V1}. Now the surface starts to expand at
first, before it shrinks. The aspect ratio of vertical to horizontal extension
is inverted during the evolution. That is, the initially tall surface becomes
flat towards the end. At the time the computations finish it appears that a
singularity may be forming soon, with the indicator $1/\mathcal{K}^m_\infty$
approaching zero faster than the surface area.

\begin{figure}
\center
\includegraphics[angle=-0,width=0.3\textwidth]{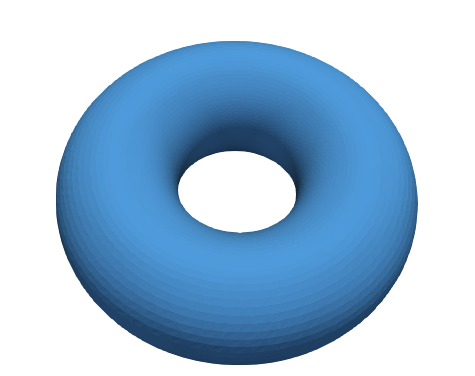} 
\includegraphics[angle=-0,width=0.3\textwidth]{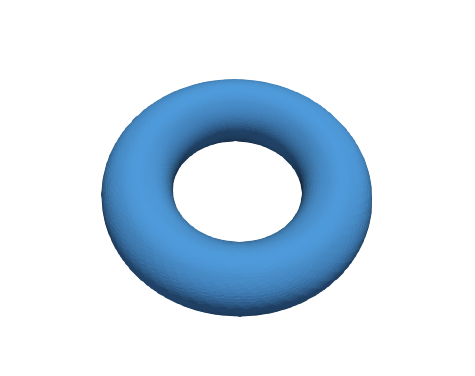} 
\includegraphics[angle=-0,width=0.3\textwidth]{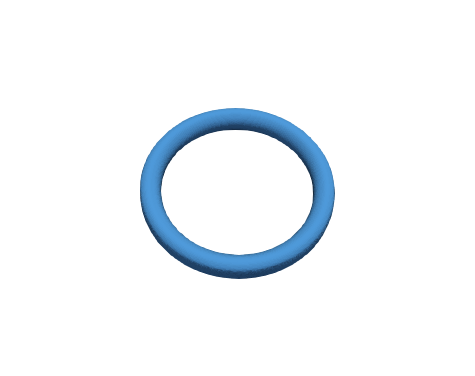} 
\includegraphics[angle=-0,width=0.3\textwidth]{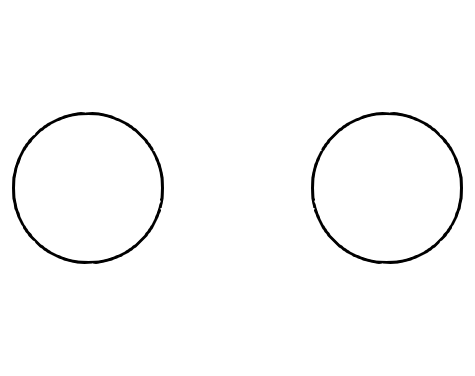} 
\includegraphics[angle=-0,width=0.3\textwidth]{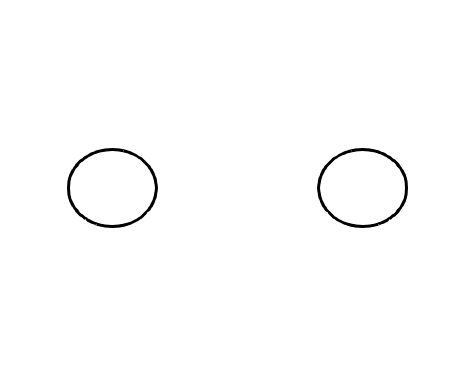} 
\includegraphics[angle=-0,width=0.3\textwidth]{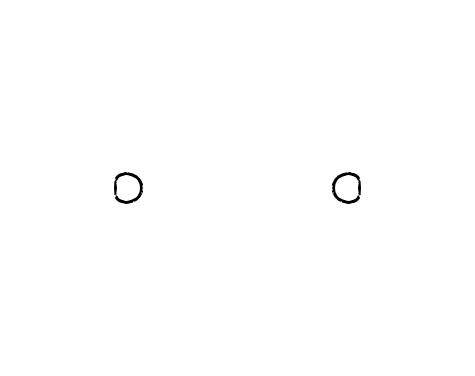} 
\\
\center
\includegraphics[angle=-90,width=0.25\textwidth]{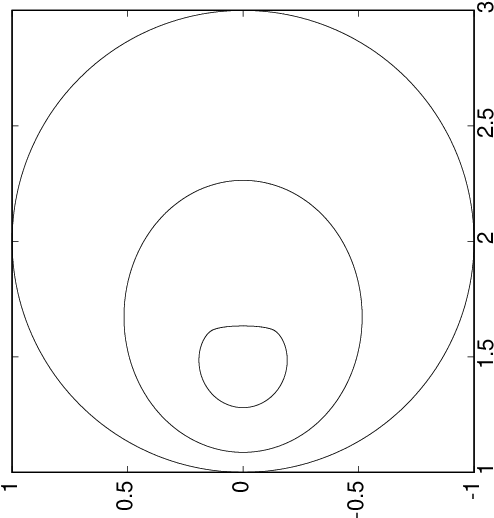}
\includegraphics[angle=-90,width=0.35\textwidth]{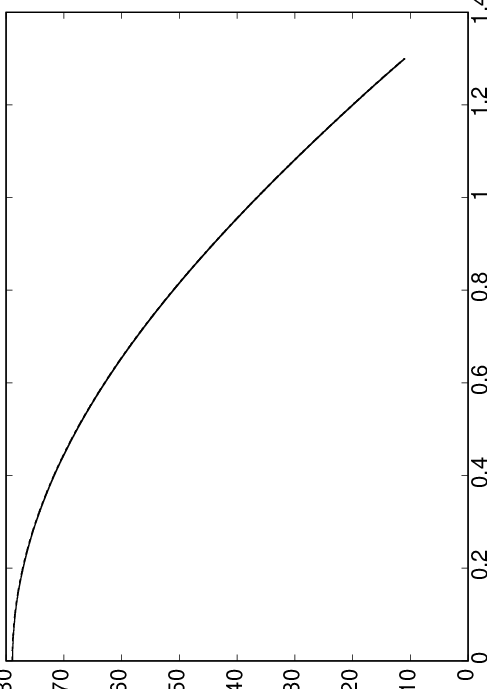}
\includegraphics[angle=-90,width=0.35\textwidth]{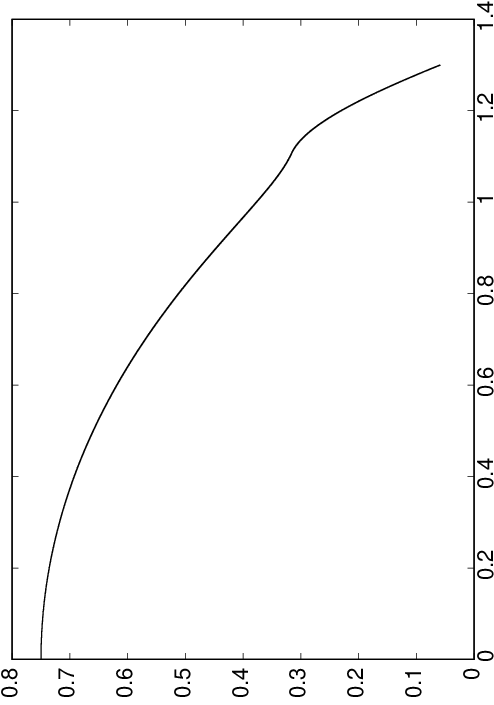}
\caption{The flow \eqref{eq:Gurtinbeta0}, with $V_0=0$, starting from a 
torus with radii $R=2$ and $r=1$.
Above we visualize $\mathcal S^m_h$ from \eqref{eq:fea_b} 
at times $t=0,1,T=1.3$ in two ways: 
a frontal view and a cut through the $x_1$-$x_2$ plane.
Below we show $x^m$ from \eqref{eq:app:fd} with $g(s)=1$ 
at times $t=0,1,T$, as well as
a plot with $|\mathcal S^m_h|$ (solid) and $A^m$ (dashed), and a plot of
$1/\mathcal{K}^m_\infty$ (right) over time.
}
\label{fig:torus21V0}
\end{figure}

A simulation for an initial torus with radii $R=2$ and $r=1$ can be seen in
Figure~\ref{fig:torus21V0} for $V_0=0$.
Here we see a singularity develop, before the torus can shrink to a circle
and vanish,
when the outer part of the torus becomes flat. 
This kind of singularity cannot occur for mean
curvature flow, and so is unique to the hyperbolic flow considered here.
In particular, we recall from \cite[Proposition~3]{SonerS93} that for mean
curvature flow it is known that there exists a critical
value $r_{crit} \in (0,1)$, such that for $0<\frac{r}{R}<r_{crit}$
the solution shrinks to a circle,
while for $r_{crit}<\frac{r}{R}<1$ it closes up the hole at the origin. 
This was also confirmed by numerical experiments in e.g.\
\cite[Figs.\ 5, 6]{gflows3d}, \cite[Figs.\ 2, 3]{aximcf} and
\cite[Figs.\ 2, 3, 4]{schemeD}, where the latter also provides an estimate for
$r_{crit}$ of about $0.6415$.

\begin{figure}
\center
\includegraphics[angle=-0,width=0.3\textwidth]{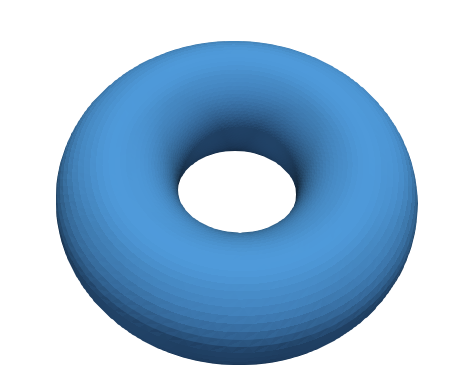}\ 
\includegraphics[angle=-0,width=0.3\textwidth]{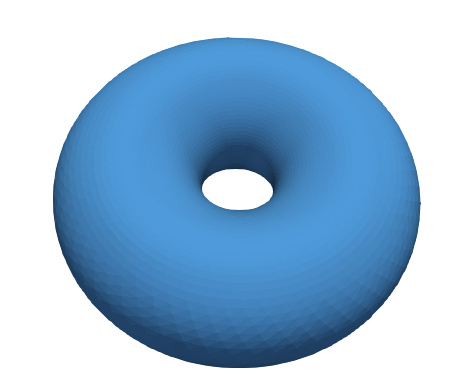}\ 
\includegraphics[angle=-0,width=0.3\textwidth]{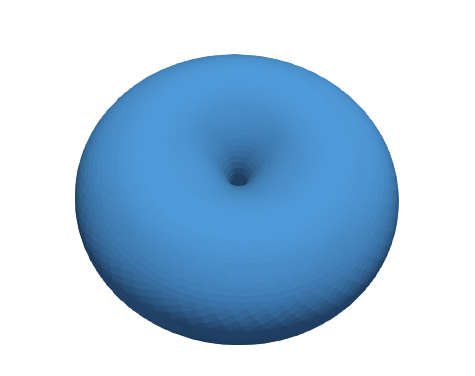} 
\includegraphics[angle=-0,width=0.3\textwidth]{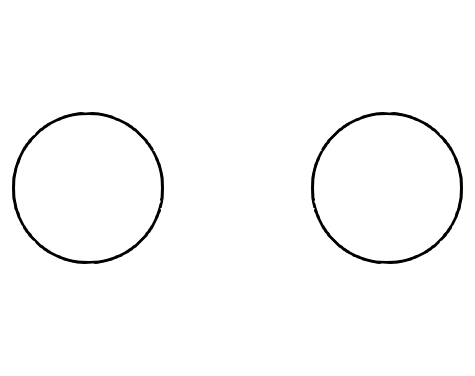}\ 
\includegraphics[angle=-0,width=0.3\textwidth]{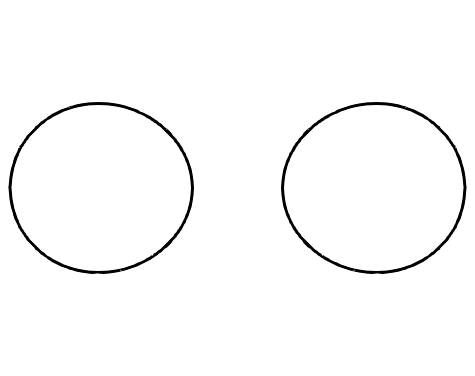}\
\includegraphics[angle=-0,width=0.3\textwidth]{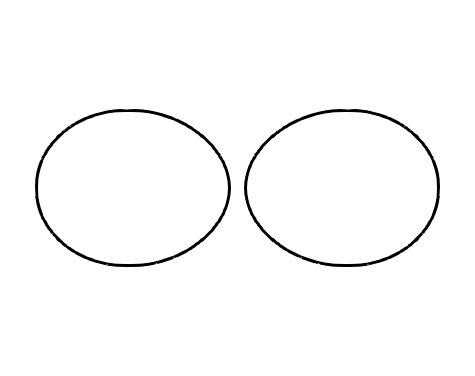} 
\\
\center
\includegraphics[angle=-90,width=0.3\textwidth]{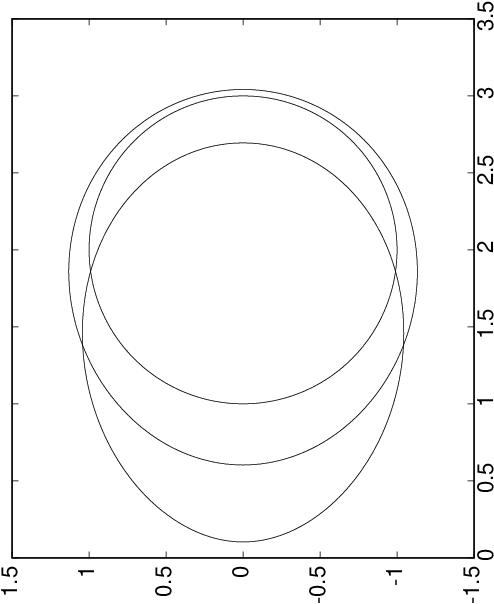}
\includegraphics[angle=-90,width=0.3\textwidth]{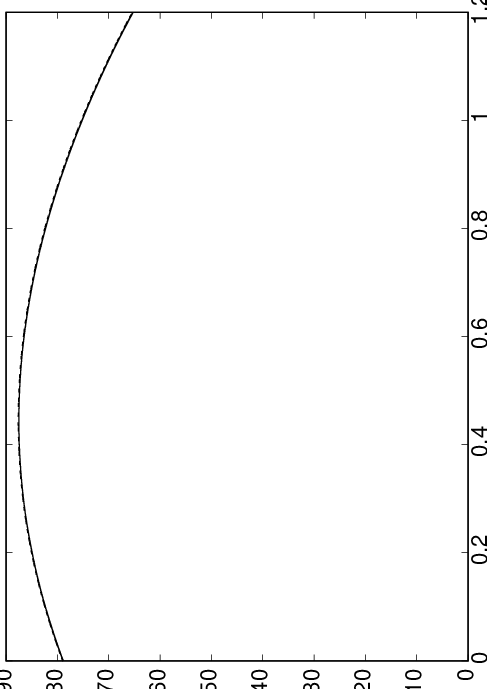}
\includegraphics[angle=-90,width=0.3\textwidth]{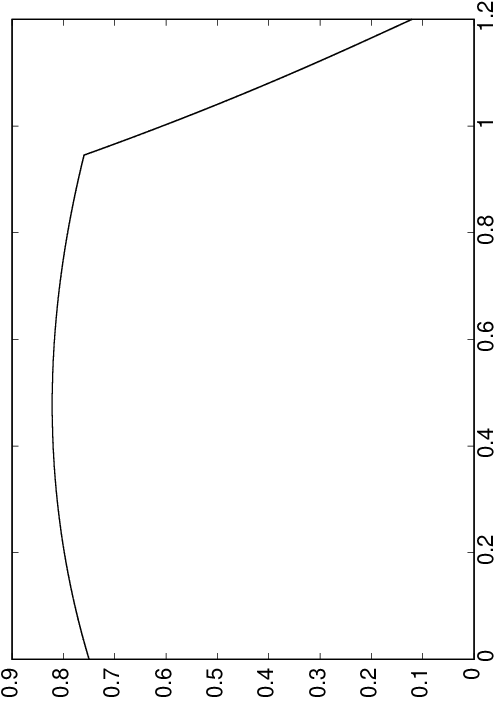}
\caption{The flow \eqref{eq:Gurtinbeta0}, with $V_0=0.5$, starting from a 
torus with radii $R=2$ and $r=1$.
Above we visualize $\mathcal S^m_h$ from \eqref{eq:fea_b} 
at times $t=0,0.7,T=1.2$ in two ways: 
a frontal view and a cut through the $x_1$-$x_2$ plane.
Below we show $x^m$ from \eqref{eq:app:fd} with $g(s)=1$ 
at times $t=0,0.7,T$, as well as
a plot with $|\mathcal S^m_h|$ (solid) and $A^m$ (dashed), and a plot of
$1/\mathcal{K}^m_\infty$ (right) over time.
}
\label{fig:torus21V05}
\end{figure}

In Figure~\ref{fig:torus21V05} we repeat the same computation but choose 
$V_0=0.5$. This results in the torus initially getting fatter, so
that in the ensuing evolution it tries to close up the inner hole to change
its topology into that of a sphere. 
This singularity is qualitatively the same as the one for mean curvature flow
when $\frac r{R} > r_{crit}$.

\begin{figure}
\center
\includegraphics[angle=0,width=0.3\textwidth]{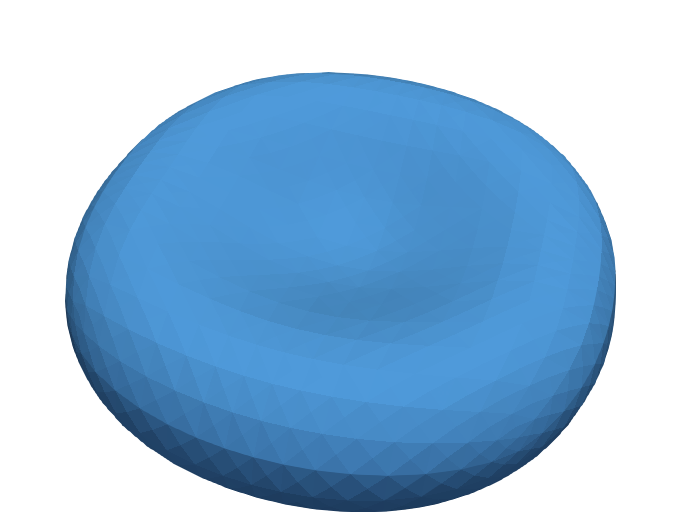} 
\includegraphics[angle=0,width=0.3\textwidth]{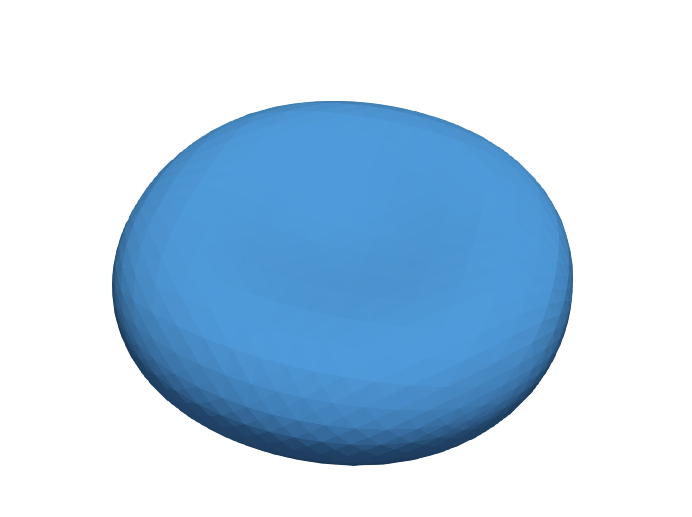} 
\includegraphics[angle=0,width=0.3\textwidth]{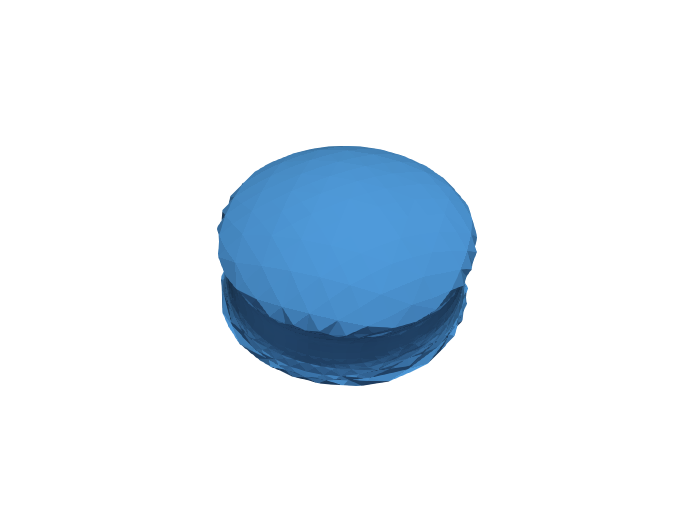} 
\includegraphics[angle=90,width=0.3\textwidth]{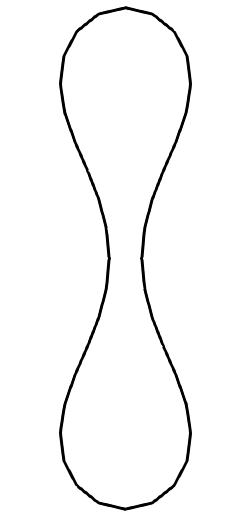} 
\includegraphics[angle=90,width=0.3\textwidth]{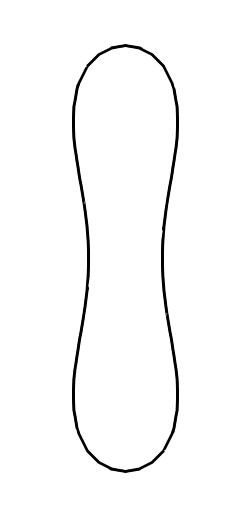} 
\includegraphics[angle=90,width=0.3\textwidth]{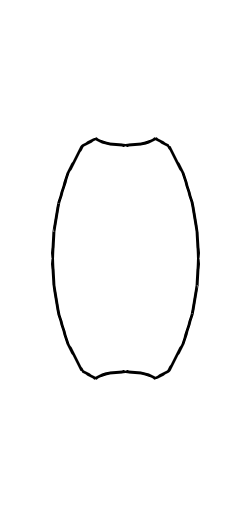} 
\caption{The flow \eqref{eq:Gurtinbeta0}, with $V_0=0$, starting from a 
biconcave disk.
We visualize $\mathcal S^m_h$ from \eqref{eq:fea_b} at times $t=0,0.5,T=0.95$
in two ways: 
a frontal view and a cut through the $x_1$-$x_2$ plane.
}
\label{fig:redbloodcellV0}
\end{figure}

We end this section with a simulation for hyperbolic mean curvature flow for a
surface that is the boundary of a nonconvex domain. 
In particular, as the initial discrete surface $\mathcal{S}^0_h$ we choose
the biconcave disk obtained in \cite[Fig.\ 8]{nsns}. 
The total dimensions of $\mathcal{S}^0_h$ are about 
$4.2 \times 4.2 \times 1.1$, and the discrete surface is made up of 3072
elements, with the largest element being about three times the size of the
smallest one.
The numerical results for
the scheme \eqref{eq:fea_b} are presented in Figure~\ref{fig:redbloodcellV0}.
We notice that at the final time $t=T$ the mesh quality has
deteriorated significantly, with the largest element now being nearly 500 times
the size of the smallest one.

\begin{figure}
\center
\includegraphics[angle=-90,width=0.5\textwidth]{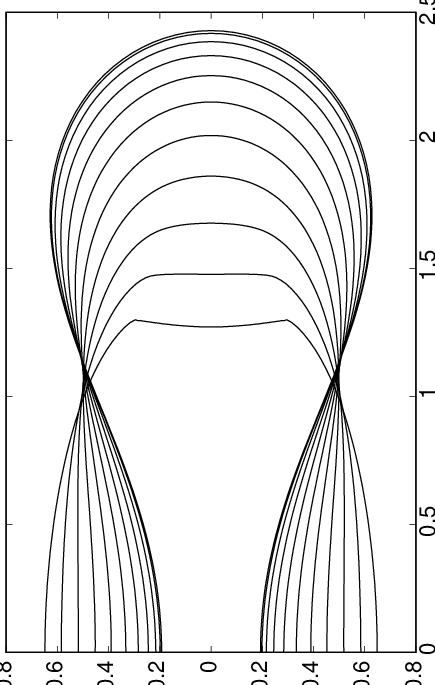}\qquad
\includegraphics[angle=-90,width=0.3\textwidth]{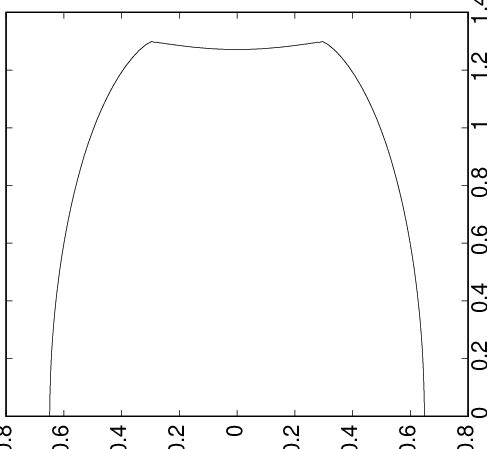} \\
\includegraphics[angle=-90,width=0.4\textwidth]{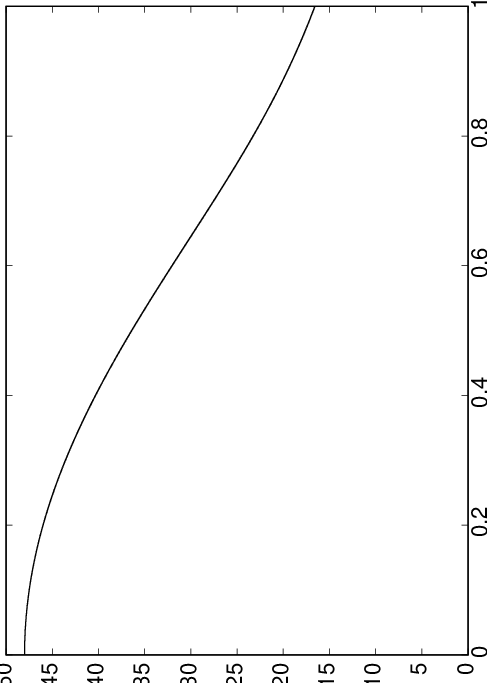}
\qquad
\includegraphics[angle=-90,width=0.4\textwidth]{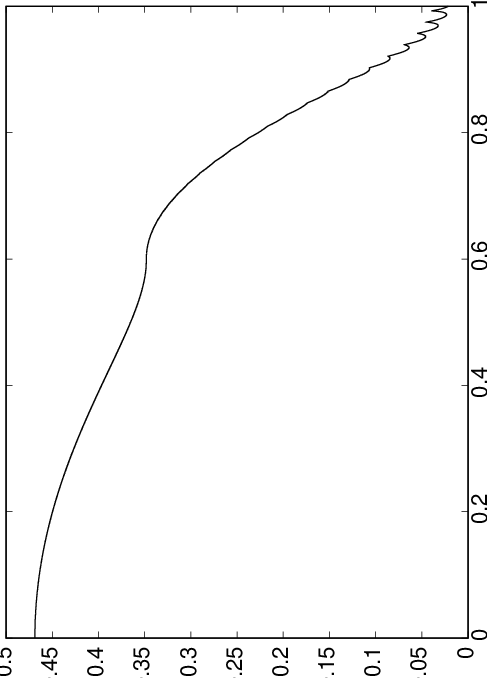}
\caption{The flow \eqref{eq:Gurtinbeta0}, with $V_0=0$, starting from a 
biconcave disk.
We show $x^m$ from \eqref{eq:app:fd} with $g(s)=1$ 
at times $t=0,0.1,\ldots,T=1$, and separately at time $t=T$.
Below we show the evolutions of $A^m$ (left) and 
$1/{\mathcal K}^m_\infty$ (right) over time.
}
\label{fig:axiredbloodcellV0}
\end{figure}

In order to qualitatively confirm the evolution shown in 
Figure~\ref{fig:redbloodcellV0}, we compute a very similar evolution for the
axisymmetric finite difference scheme \eqref{eq:app:fd} with $g(s)=1$. Here we choose
as initial data $x^0$ the curve displayed in the middle of
\cite[Fig.\ 5]{axipwf}.
The obtained simulation is shown in Figure~\ref{fig:axiredbloodcellV0},
which agrees very well with the full 3d results in 
Figure~\ref{fig:redbloodcellV0}.
In particular, at the side of the disk a concavity appears, which then leads to
a blow-up in curvature and a resulting singularity for the flow.

\subsection{The hyperbolic mean curvature flow \eqref{eq:LeFloch}}

\subsubsection{Radial evolutions for a sphere}

We begin with a convergence experiment in order to validate our schemes 
\eqref{eq:feaLF} and \eqref{eq:app:fd} with $g(s)=1+\tfrac12 s$, using the solution discussed in
\cite[\S7.1]{LeFlochS08}. 
In fact,  a family of spheres in $\bR^3$ with 
radius $r(t)$ and outer normal $\nu(t)$ is a solution to \eqref{eq:LeFloch} 
with $\mathcal V_{\mathcal{S} \mid_{t=0}} = V_0 \in \bR$ if
\begin{equation} \label{eq:rtt}
\ddot{r}(t) = - \frac{1}{r(t)} \bigl( (\dot{r}(t))^2 +2 \bigr), \quad \text{in } (0,T],\quad
r(0) = r_0>0,\ \dot{r}(0) = V_0.
\end{equation}
It can be shown that the solution of \eqref{eq:rtt} is given by
\begin{equation} \label{eq:rtrue}
r(t) = \sqrt{r_0^2 + 2r_0V_0t - 2t^2}\,.
\end{equation}
Thus, if $V_0 > 0$, then $r(t)$ will at first increase until it hits
a maximum at $t=\tfrac12 r_0V_0$ after which it  
will decrease and shrink to a point at time $T_{max}$, where $T_{max}=\frac{r_0}{2} \bigl( V_0 + \sqrt{V_0^2 +2} \bigr)$, see
\cite[p.~611]{LeFlochS08}. 
On the other hand, 
if $V_0 \leq 0$ then the sphere will monotonically shrink to a point at time $T_{max}$. 

We use the true solution \eqref{eq:rtrue} 
for a convergence test to validate the 
schemes \eqref{eq:feaLF} and \eqref{eq:app:fd} with $g(s)=1+\tfrac12 s$. 
We begin with the scheme \eqref{eq:feaLF}, using
a sequence of triangulations of the unit sphere.  The results
for $V_0 = 0$, $V_0 = 1$ and $V_0 = -1$ are reported in Table~\ref{tab:3dLF}. 
Here we employed the time step size choice $\Delta t = 0.25 h_0$.
In contrast to Table~\ref{tab:3d}, the observed order of convergence
is no longer quadratic, which is likely due to the explicit approximation of 
$g(\mathcal V_{\mathcal S}^2)$ in \eqref{eq:feaLF}.
\begin{table}
\center
\begin{tabular}{|r|c|c|c|c|c|c|c|}
\hline 
& & \multicolumn{2}{c|}{$V_0 = 0$} & \multicolumn{2}{c|}{$V_0 = 1$} 
& \multicolumn{2}{c|}{$V_0 = -1$} \\ 
$J$   & $h_0$     & $\errorSS$&EOC & $\errorSS$&EOC & $\errorSS$&EOC  \\ \hline
1536  & 2.0854e-01& 4.7827e-03&1.66& 1.1273e-02&1.89& 2.5920e-02&--- \\ 
6144  & 1.0472e-01& 1.3429e-03&1.84& 4.4218e-03&1.36& 1.2374e-02&1.07\\ 
24576 & 5.2416e-02& 4.0775e-04&1.72& 1.4151e-03&1.65& 6.0207e-03&1.04\\ 
98304 & 2.6215e-02& 1.3576e-04&1.59& 5.7473e-04&1.30& 2.6327e-03&1.19\\ 
393216& 1.3108e-02& 5.7032e-05&1.25& 2.4144e-04&1.25& 1.2466e-03&1.08\\ 
\hline
\end{tabular}
\caption{Errors for \eqref{eq:feaLF} for the convergence test for 
\eqref{eq:rtrue} with $r_0 = 1$ and $V_0 = 0$, $V_0 = 1$ or $V_0 = -1$,
over the time interval $[0,0.25]$.
We also display the experimental orders of convergence (EOC).
}
\label{tab:3dLF}
\end{table}%

Moreover, we perform the same convergence experiment also for the axisymmetric
scheme \eqref{eq:app:fd} with $g(s)=1+\tfrac12 s$. To this end, we choose as initial data
$x_0(\rho) = (\cos (\frac32\pi + \pi \rho), \sin (\frac32\pi + \pi \rho) )^T$. Similarly as above,
the mapping $x(\rho,t)=r(t) x_0(\rho)$ is an exact solution
of \eqref{eq:app:hcsf}.
The observed errors \eqref{eq:errorXX} for the choice $\Delta t =h$ are shown in Table~\ref{tab:axiLF}.
Once again we note that in contrast to Table~\ref{tab:axi}, the observed
convergence rates now are only linear. This is very likely due to the
effect of the coefficient function $g$ in \eqref{eq:app:fda} for the choice
$g(s)=1+\tfrac12 s$, compared to the constant $g(s) = 1$.
\begin{table}
\center
\begin{tabular}{|r|c|c|c|c|c|c|}
\hline
& \multicolumn{2}{c|}{$V_0 = 0$} & \multicolumn{2}{c|}{$V_0 = 1$} 
& \multicolumn{2}{c|}{$V_0 = -1$} \\ 
$J$ & $\errorXX$&EOC & $\errorXX$&EOC & $\errorXX$&EOC  \\ \hline
32  & 5.9402e-03&---  &5.6974e-03 &---  &1.0942e-02&---  \\
64  & 2.8216e-03&1.07 &2.5537e-03 &1.16 &5.1825e-03&1.08 \\
128 & 1.3817e-03&1.03 &1.2096e-03 &1.08 &2.5339e-03&1.03 \\ 
256 & 6.8508e-04&1.01 &5.8875e-04 &1.04 &1.2551e-03&1.01 \\ 
512 & 3.4135e-04&1.01 &2.9047e-04 &1.02 &6.2499e-04&1.01 \\ 
\hline
\end{tabular}
\caption{
Errors for \eqref{eq:app:fd} with $g(s)=1+\tfrac12 s$ for the convergence test for 
\eqref{eq:rtrue} with $r_0=1$ and $V_0=0$, $V_0 = 1$ or $V_0 = 1$,
over the time intervals $[0,0.5]$, $[0,0.5]$ and $[0,0.25]$, respectively.
We also display the experimental orders of convergence (EOC).
}
\label{tab:axiLF}
\end{table}%
To further investigate the suboptimal convergence rates, we repeat the
experiments for Tables~\ref{tab:3dLF} and \ref{tab:axiLF} but now choose 
$\Delta t = 0.5 h^2_0$ and $\Delta t = h^2$, respectively.
The results are shown in Tables~\ref{tab:3dLFp2} and \ref{tab:axiLFp2},
demonstrating quadratic convergence rates.
Hence, as expected, these results appear to confirm that the
temporal discretization, and in particular the explicit approximation of 
$g(\mathcal V_{\mathcal S}^2)$, is responsible for the reduced convergence
rates seen in Tables~\ref{tab:3dLF} and \ref{tab:axiLF}. 
\begin{table}
\center
\begin{tabular}{|r|c|c|c|c|c|c|c|}
\hline 
& & \multicolumn{2}{c|}{$V_0 = 0$} & \multicolumn{2}{c|}{$V_0 = 1$} 
& \multicolumn{2}{c|}{$V_0 = -1$} \\ 
$J$   & $h_0$     & $\errorSS$&EOC & $\errorSS$&EOC & $\errorSS$&EOC  \\ \hline
1536  & 2.0854e-01& 4.6490e-03&1.62& 7.2549e-03&2.12& 1.4806e-02&--- \\ 
6144  & 1.0472e-01& 1.5653e-03&1.58& 2.1150e-03&1.79& 3.3757e-03&2.15\\ 
24576 & 5.2416e-02& 4.4500e-04&1.82& 6.3384e-04&1.74& 8.1195e-04&2.06\\ 
98304 & 2.6215e-02& 1.2257e-04&1.86& 1.8075e-04&1.81& 2.2052e-04&1.88\\ 
393216& 1.3108e-02& 3.3678e-05&1.86& 5.0157e-05&1.85& 5.9463e-05&1.89\\ 
\hline
\end{tabular}
\caption{Same as Table~\ref{tab:3dLF} but with $\Delta t= 0.5 h^2_0$.}
\label{tab:3dLFp2}
\end{table}%
\begin{table}
\center
\begin{tabular}{|r|c|c|c|c|c|c|}
\hline
& \multicolumn{2}{c|}{$V_0 = 0$} & \multicolumn{2}{c|}{$V_0 = 1$} 
& \multicolumn{2}{c|}{$V_0 = -1$} \\ 
$J$ & $\errorXX$&EOC & $\errorXX$&EOC & $\errorXX$&EOC  \\ \hline
32  &4.1126e-04 &--- &5.8329e-05 &---  &4.4130e-04 &--- \\
64  &1.0181e-04 &2.01&6.9665e-06 &3.07 &1.1113e-04 &1.99\\
128 &2.5403e-05 &2.00&1.1422e-06 &2.61 &2.6938e-05 &2.04\\ 
256 &6.3444e-06 &2.00&2.1542e-07 &2.41 &6.7589e-06 &1.99\\ 
512 &1.5813e-06 &2.00&4.2844e-08 &2.33 &1.6773e-06 &2.01\\ 
\hline
\end{tabular}
\caption{Same as Table~\ref{tab:axiLF} but with $\Delta t= h^2$.}
\label{tab:axiLFp2}
\end{table}%

\begin{figure}
\center
\includegraphics[angle=-90,width=0.3\textwidth]{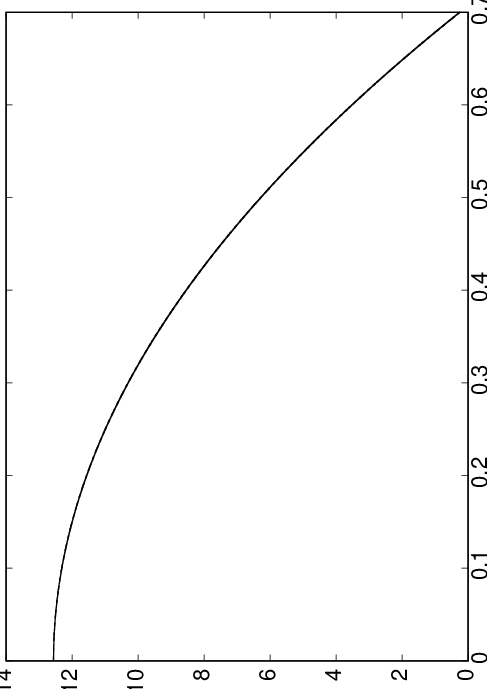}
\includegraphics[angle=-90,width=0.3\textwidth]{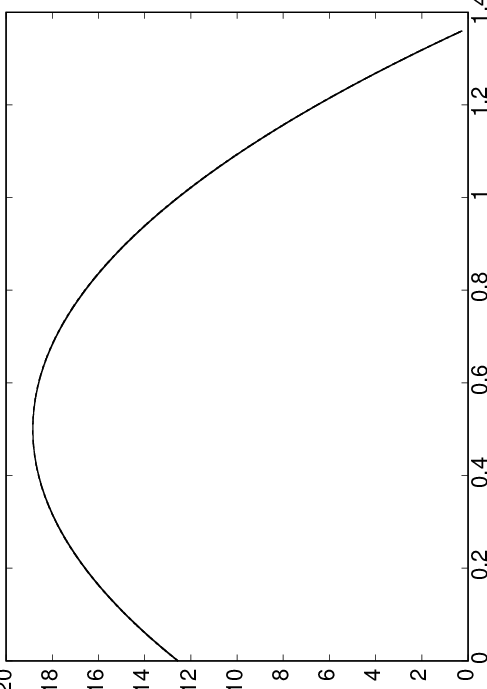}
\includegraphics[angle=-90,width=0.3\textwidth]{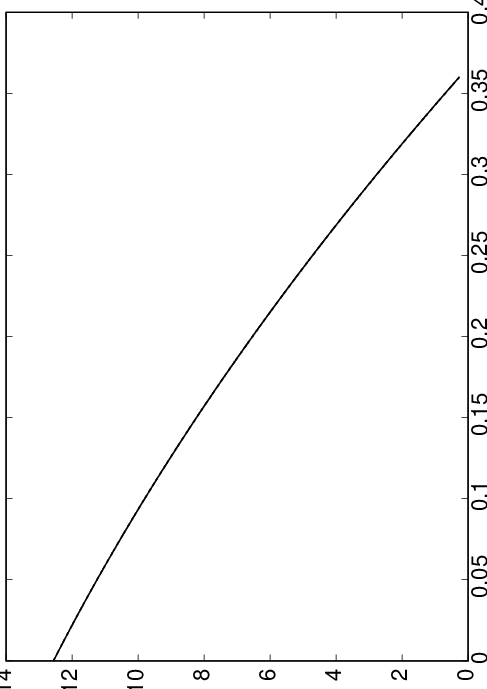}
\\
\includegraphics[angle=-90,width=0.3\textwidth]{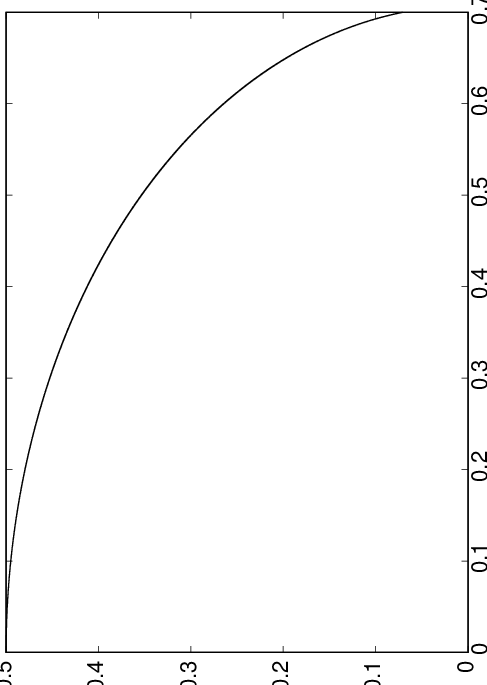}
\includegraphics[angle=-90,width=0.3\textwidth]{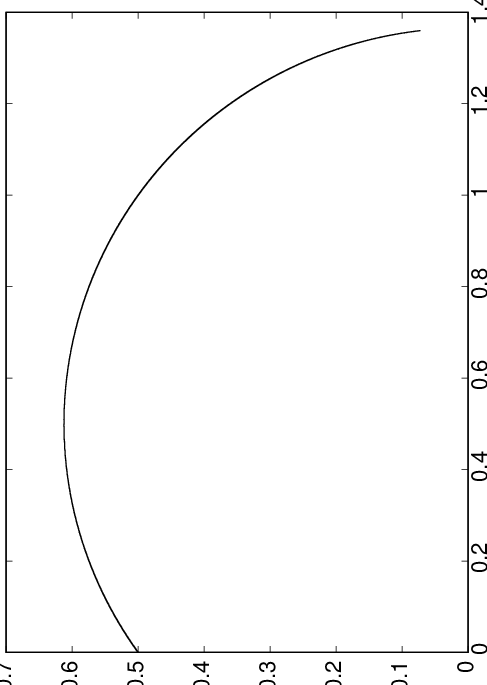}
\includegraphics[angle=-90,width=0.3\textwidth]{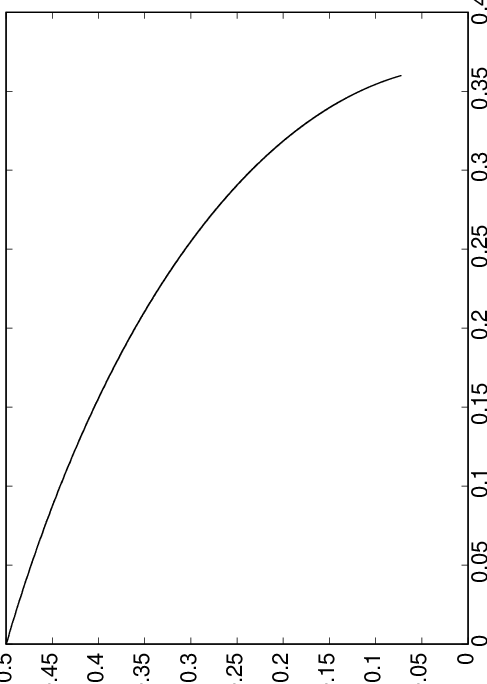}
\caption{The flow \eqref{eq:LeFloch} starting from a unit sphere.
Above we show the evolutions of $|\mathcal S^m_h|$ (solid) and $A^m$ (dashed) 
over time for $V_0=0$, $V_0=1$ and $V_0=-1$.
The final times for these computations are $T=0.7$, $T=1.36$ and $T=0.36$,
respectively.
Below we show the corresponding evolutions of $1/{\mathcal K}^m_\infty$ over 
time.
}
\label{fig:axispheresLF}
\end{figure}

As discussed above, for an initial sphere, depending on the sign of 
$V_0 \in \bR$, the family of
spheres either expands at first and then shrinks, or shrinks immediately. We
visualize these different behaviours in Figure~\ref{fig:axispheresLF}. 
In each case we observe a smooth solution until the spheres shrink to a point, 
meaning that $|\mathcal S^m_h|$ and $1/{\mathcal K}^m_\infty$ approach 
zero at the same time.
Compared to Figure~\ref{fig:axispheres} we note small differences in 
the time scales,
and in the maximal radius reached for the simulation with $V_0=1$. 

\subsubsection{Axisymmetric evolutions for nonspherical surfaces}

In this subsection we repeat some of the experiments in \S\ref{sec:nraxi} but 
now for the flow \eqref{eq:LeFloch}, 
using the schemes \eqref{eq:feaLF} and \eqref{eq:app:fd} with $g(s)=1+\tfrac12 s$.

\begin{figure}
\center
\includegraphics[angle=-90,width=0.15\textwidth]{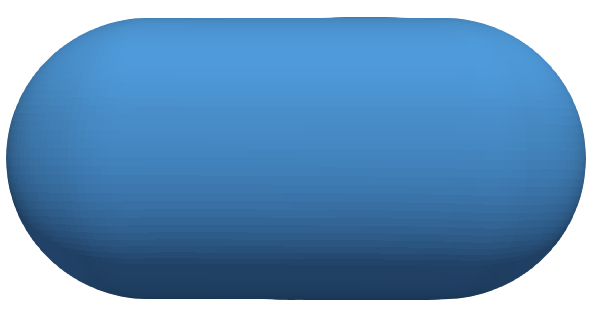} 
\includegraphics[angle=-90,width=0.15\textwidth]{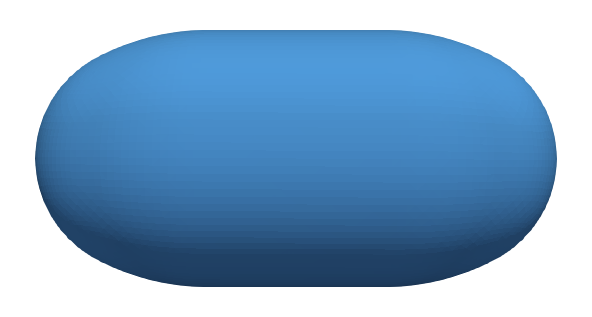} 
\includegraphics[angle=-90,width=0.15\textwidth]{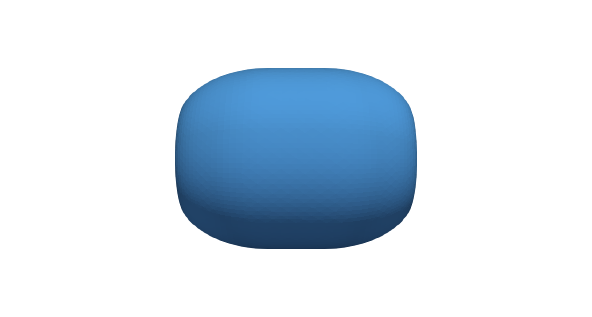} 
\includegraphics[angle=-90,width=0.15\textwidth]{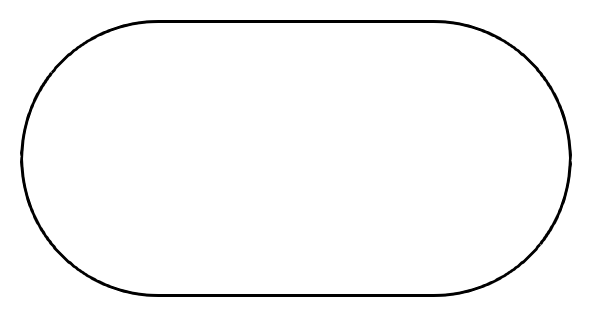} 
\includegraphics[angle=-90,width=0.15\textwidth]{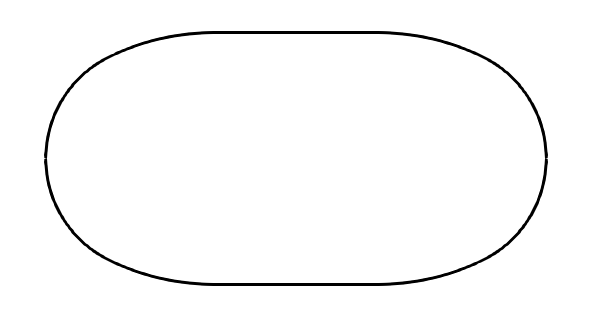} 
\includegraphics[angle=-90,width=0.15\textwidth]{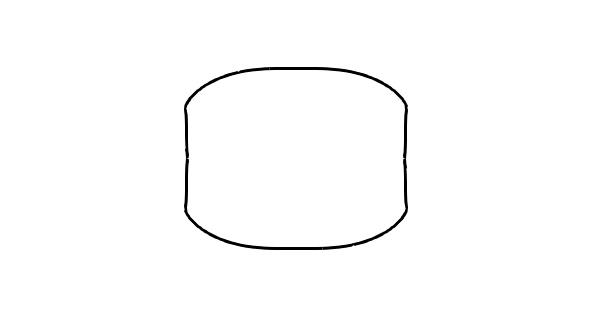} 
\center
\includegraphics[angle=-90,width=0.12\textwidth]{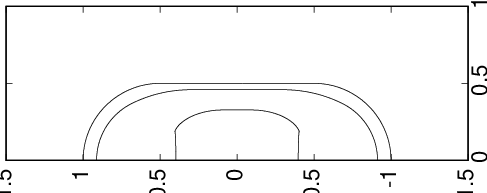}
\includegraphics[angle=-90,width=0.4\textwidth]{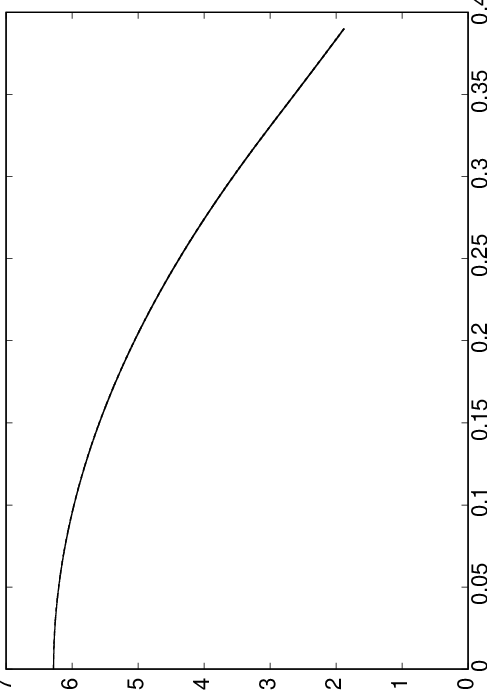}
\includegraphics[angle=-90,width=0.4\textwidth]{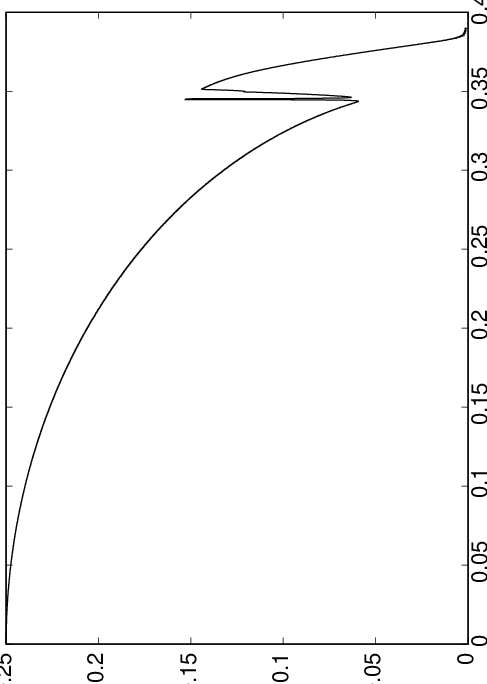}
\caption{The flow \eqref{eq:LeFloch}, with $V_0=0$, starting from a 
$2:1:1$ cigar.
Above we visualize $\mathcal S^m_h$ from \eqref{eq:feaLF} 
at times $t=0,0.2,T=0.39$ in two ways: a frontal view and a cut through the
$x_1$-$x_2$ plane.
Below we show $x^m$ from \eqref{eq:app:fd} with $g(s)=1+\tfrac12 s$ 
at times $t=0,0.2,T$, as well as
a plot with $|\mathcal S^m_h|$ (solid) and $A^m$ (dashed), and a plot of
$1/\mathcal{K}^m_\infty$ (right) over time.
}
\label{fig:cigar211LFV0}
\end{figure}

In analogy to the simulation in Figure~\ref{fig:cigar211V0}, we
begin with an experiment for an initial $2:1:1$ cigar,
when the initial velocity is chosen as $V_0=0$,
see Figure~\ref{fig:cigar211LFV0}. 
Interestingly, the ensuing evolution is slightly different. 
Apart from the shorter time scale, the final shape of the surface
appears to be nonconvex.

\begin{figure}
\center
\includegraphics[angle=-0,width=0.3\textwidth]{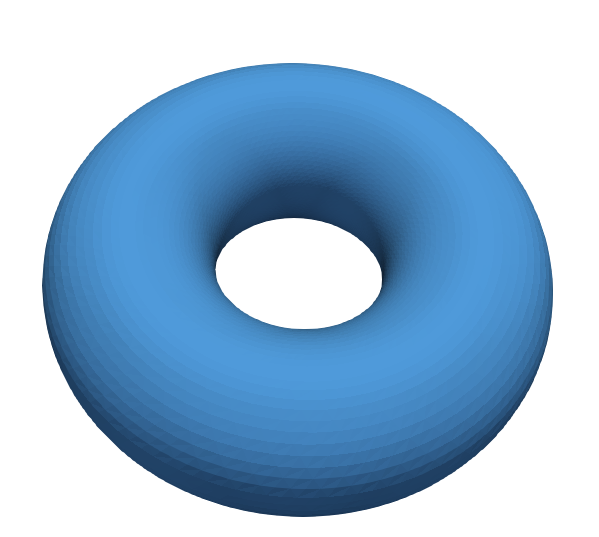} 
\includegraphics[angle=-0,width=0.3\textwidth]{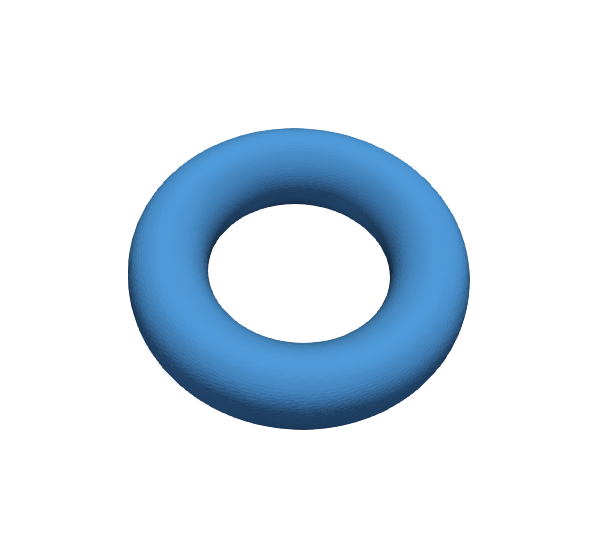} 
\includegraphics[angle=-0,width=0.3\textwidth]{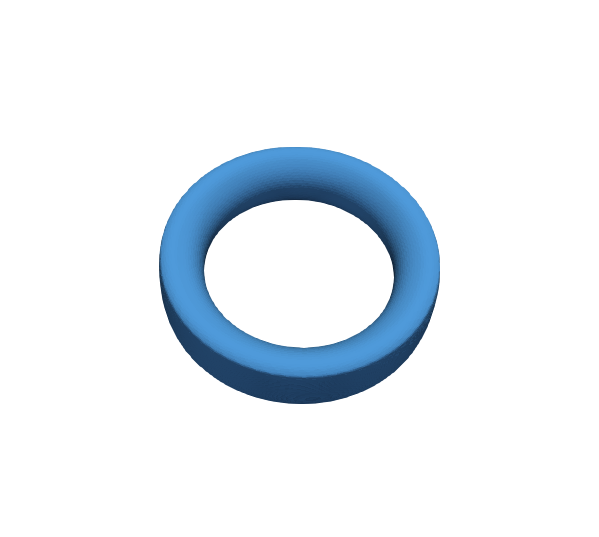} 
\includegraphics[angle=-0,width=0.3\textwidth]{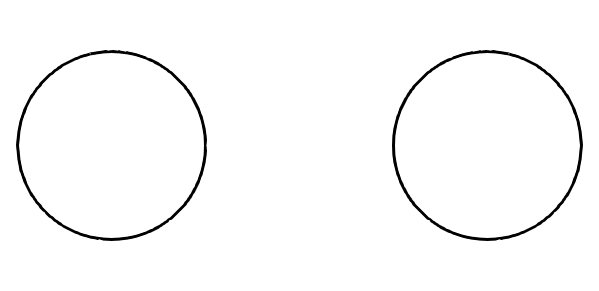} 
\includegraphics[angle=-0,width=0.3\textwidth]{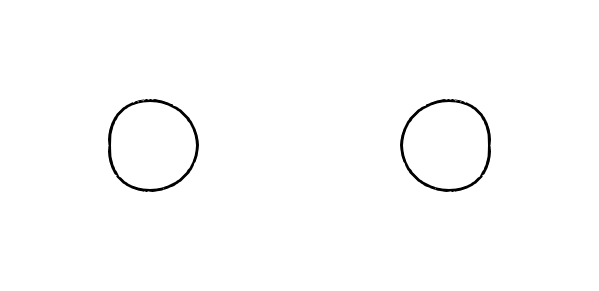} 
\includegraphics[angle=-0,width=0.3\textwidth]{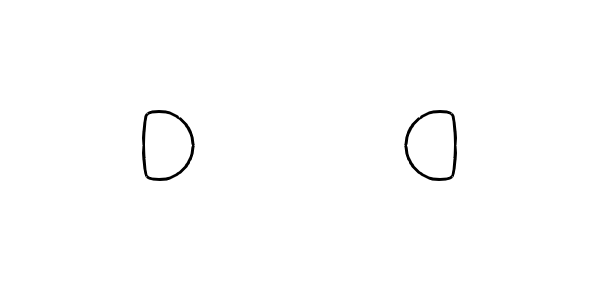} 
\\
\center
\includegraphics[angle=-90,width=0.25\textwidth]{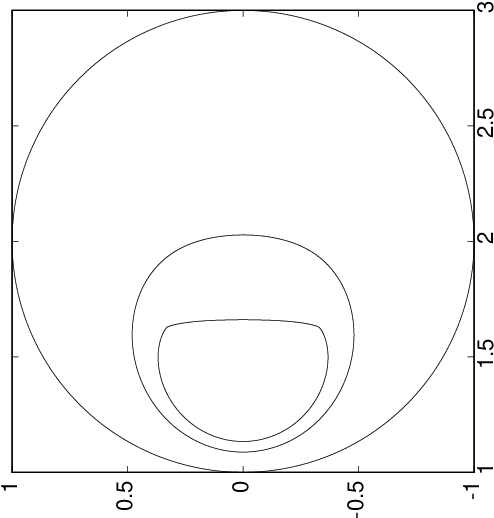}
\includegraphics[angle=-90,width=0.35\textwidth]{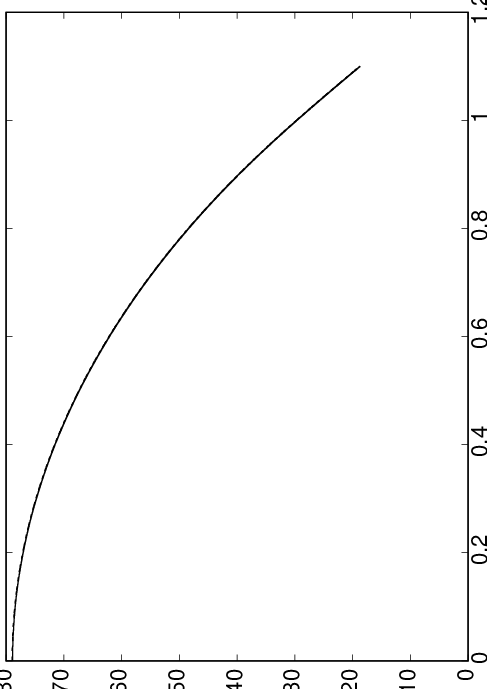}
\includegraphics[angle=-90,width=0.35\textwidth]{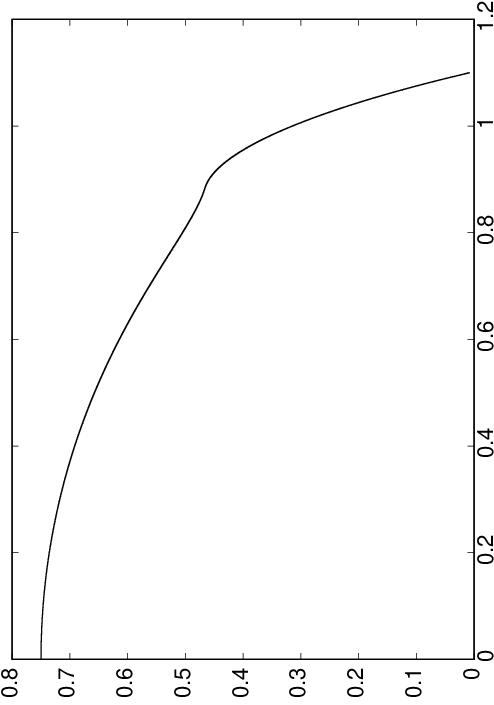}
\caption{The flow \eqref{eq:LeFloch}, with $V_0=0$, starting from a 
torus with radii $R=2$ and $r=1$.
Above we visualize $\mathcal S^m_h$ from \eqref{eq:feaLF} 
at times $t=0,1,T=1.1$ in two ways: 
a frontal view and a cut through the $x_1$-$x_2$ plane.
Below we show $x^m$ from \eqref{eq:app:fd} with $g(s)=1+\tfrac12 s$ 
at times $t=0,1,T$, as well as
a plot with $|\mathcal S^m_h|$ (solid) and $A^m$ (dashed), and a plot of
$1/\mathcal{K}^m_\infty$ (right) over time.
}
\label{fig:torus21LFV0}
\end{figure}

\begin{figure}
\center
\includegraphics[angle=-0,width=0.3\textwidth]{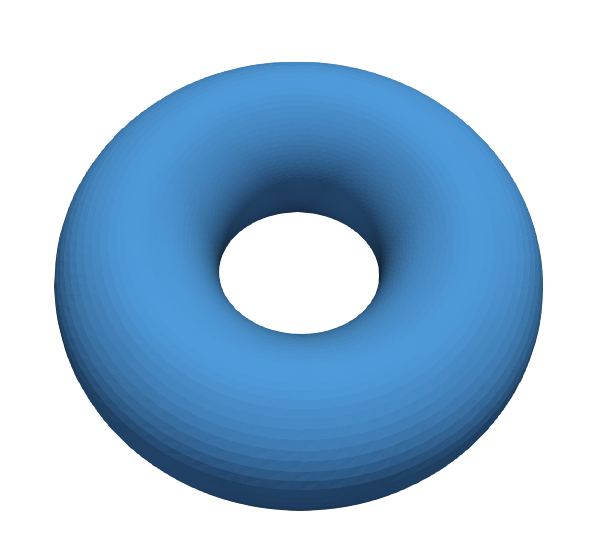}\ 
\includegraphics[angle=-0,width=0.3\textwidth]{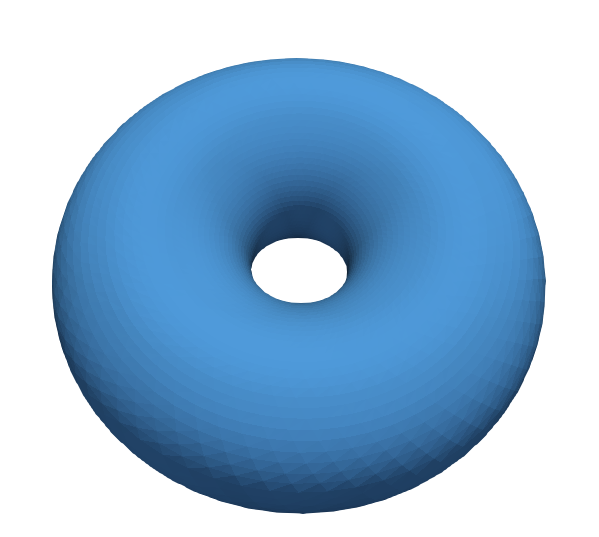}\ 
\includegraphics[angle=-0,width=0.3\textwidth]{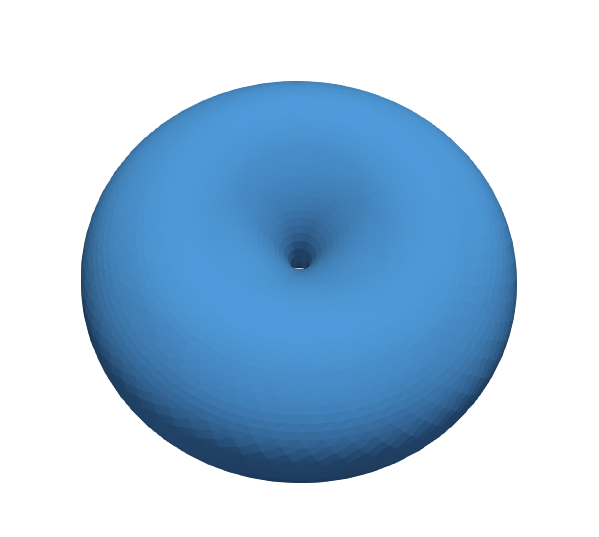} 
\includegraphics[angle=-0,width=0.3\textwidth]{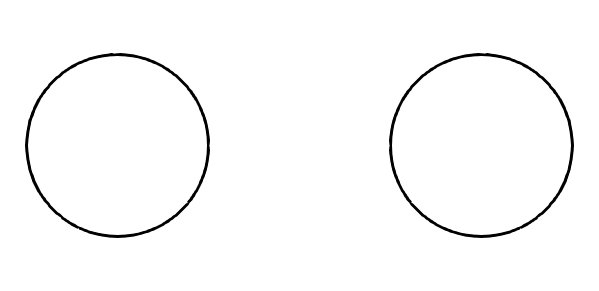}\ 
\includegraphics[angle=-0,width=0.3\textwidth]{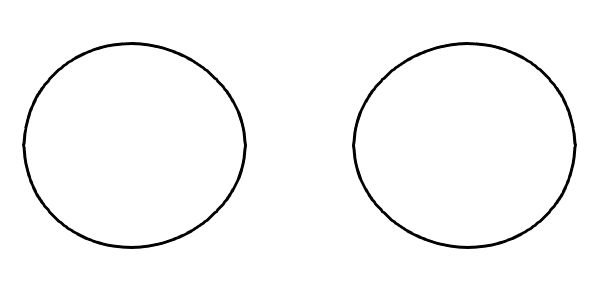}\
\includegraphics[angle=-0,width=0.3\textwidth]{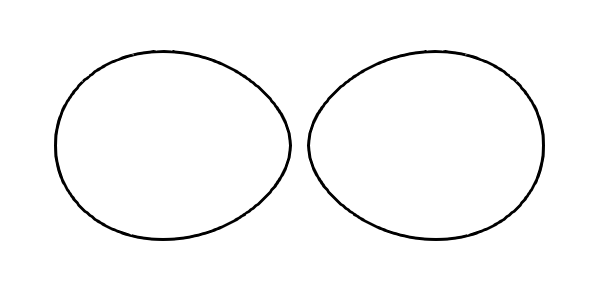} 
\\
\center
\includegraphics[angle=-90,width=0.3\textwidth]{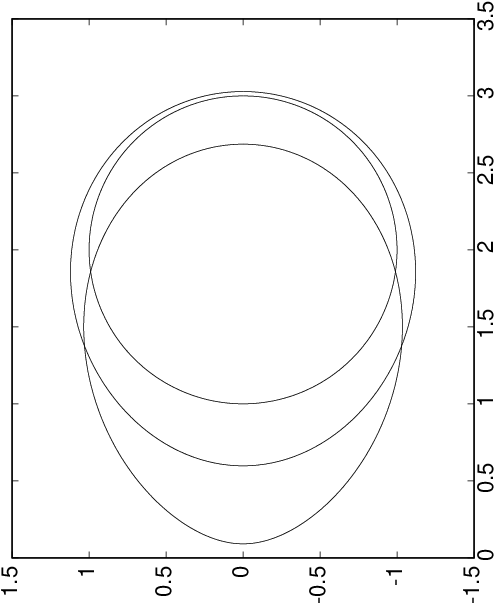}
\includegraphics[angle=-90,width=0.3\textwidth]{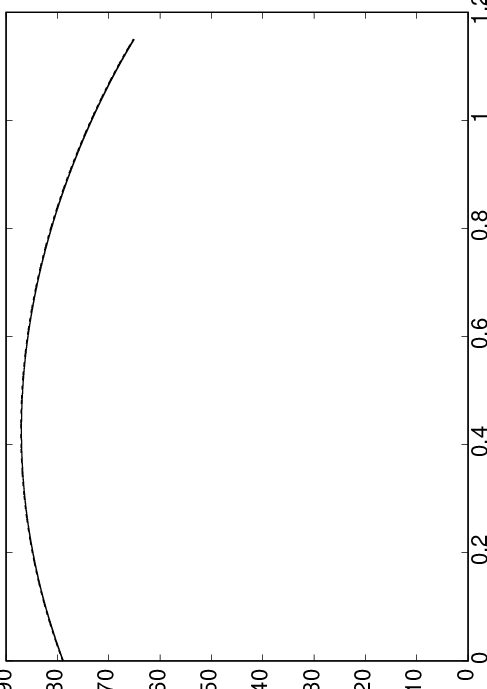}
\includegraphics[angle=-90,width=0.3\textwidth]{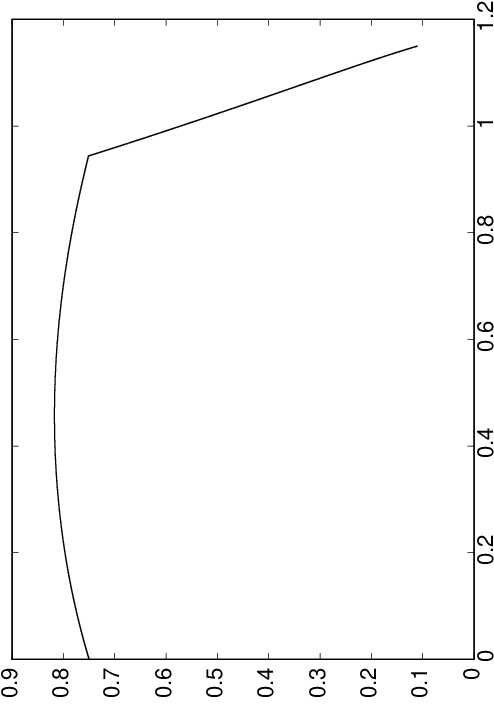}
\caption{The flow \eqref{eq:LeFloch}, with $V_0=0.5$, starting from a 
torus with radii $R=2$ and $r=1$.
Above we visualize $\mathcal S^m_h$ from \eqref{eq:feaLF} 
at times $t=0,0.5,T=1.15$ in two ways: 
a frontal view and a cut through the $x_1$-$x_2$ plane.
Below we show $x^m$ from \eqref{eq:app:fd} with $g(s)=1+\tfrac12 s$ 
at times $t=0,0.5,T$, as well as
a plot with $|\mathcal S^m_h|$ (solid) and $A^m$ (dashed), and a plot of
$1/\mathcal{K}^m_\infty$ (right) over time.
}
\label{fig:torus21LFV05}
\end{figure}

We end this section with the two experiments for an initial torus with radii 
$R=2$ and $r=1$. Overall the evolutions shown in
Figures~\ref{fig:torus21LFV0} and \ref{fig:torus21LFV05} are very similar to
the earlier ones for the flow \eqref{eq:Gurtinbeta0}, recall
Figures~\ref{fig:torus21V0} and \ref{fig:torus21V05}. However, there are minor
differences in the time scales, and the shape of the surface as it approaches
the singularity for the experiment with $V_0$ is noticeably different.

\section{Conclusions}
We have presented fully practical finite element and finite difference 
approximations for two hyperbolic mean curvature flows of hypersurfaces
in $\bR^3$ that are of interest
in application and in differential geometry. 
To our knowledge, the proposed algorithms for
the flows \eqref{eq:Gurtinbeta0} and \eqref{eq:LeFloch} 
are the first numerical approximations of these flows in the literature.
Our presented numerical results confirm the onset of nontrivial singularities,
as they have been analytically predicted for analogous evolution laws for
curves in the plane.

\def\soft#1{\leavevmode\setbox0=\hbox{h}\dimen7=\ht0\advance \dimen7
  by-1ex\relax\if t#1\relax\rlap{\raise.6\dimen7
  \hbox{\kern.3ex\char'47}}#1\relax\else\if T#1\relax
  \rlap{\raise.5\dimen7\hbox{\kern1.3ex\char'47}}#1\relax \else\if
  d#1\relax\rlap{\raise.5\dimen7\hbox{\kern.9ex \char'47}}#1\relax\else\if
  D#1\relax\rlap{\raise.5\dimen7 \hbox{\kern1.4ex\char'47}}#1\relax\else\if
  l#1\relax \rlap{\raise.5\dimen7\hbox{\kern.4ex\char'47}}#1\relax \else\if
  L#1\relax\rlap{\raise.5\dimen7\hbox{\kern.7ex
  \char'47}}#1\relax\else\message{accent \string\soft \space #1 not
  defined!}#1\relax\fi\fi\fi\fi\fi\fi}

\end{document}